\documentclass{elsarticle}[preprint,review,10pt]

\journal{Journal of Computational Physics}

\usepackage{tikz}
\usetikzlibrary{arrows.meta, positioning}
\usepackage{graphicx} 

\usepackage{graphicx} 
\usepackage{algorithm}
\usepackage{algorithmic}
\usepackage{booktabs} 
\usepackage{siunitx} 
\usepackage{caption} 
\usepackage{multirow} 
\usepackage{threeparttable}
\usepackage{tabularx}
\usepackage{array}
\usepackage{adjustbox}

\usepackage[utf8]{inputenc}
\usepackage{graphicx}
\usepackage{bm}
\usepackage{comment}
\usepackage{psfrag}
\usepackage{latexsym,amsmath,amsfonts,amscd,amsthm}
\usepackage{changebar}
\usepackage{color}
\usepackage{bm}
\usepackage{tikz}
\usepackage{multirow}
\usepackage{xcolor}
\usepackage{hyperref}
\usepackage{setspace}
\usepackage{amssymb}
\usepackage{mathrsfs}
\usepackage{caption}
\usepackage{subcaption}

\usepackage[title]{appendix}
\usepackage{ulem}
\usepackage{tikz}
\usetikzlibrary{arrows,backgrounds}

\newlength{\spse}
\newtheorem{thm}{Theorem}[section]

\newtheorem{rem}[thm]{Remark}

\newcommand{\BU}{\mathbf{U}}

\newcommand{\BK}{\mathbf{K}}
\newcommand{\BI}{\mathbf{I}}

\newcommand{\BF}{\mathbf{F}}

\newcommand{\BA}{\mathbf{A}}
\newcommand{\WBA}{\widetilde{\mathbf{A}}}
\newcommand{\BQ}{\mathbf{Q}}

\newcommand{\Bupsilon}{\boldsymbol{\upsilon}}

\newcommand{\bw}{\mathbf{w}}
\newcommand{\bc}{\mathbf{c}}
\newcommand{\bb}{\mathbf{b}}
\newcommand{\wbb}{\widetilde{\mathbf{b}}}

\newcommand{\bn}{\boldsymbol{n}}

\newcommand{\brho}{\boldsymbol{\rho}}
\newcommand{\bg}{\boldsymbol{g}}
\newcommand{\br}{\boldsymbol{r}}
\newcommand{\bq}{\boldsymbol{q}}

\newcommand{\bQ}{\boldsymbol{Q}}

\newcommand{\bff}{\boldsymbol{f}}

\newcommand{\beeta}{\boldsymbol{\eta}}

\newcommand{\BC}{{\bf{C}} }
\newcommand{\BW}{{\bf{W}} }
\newcommand{\BM}{{\bf{M}} }

\newcommand{\BD}{{\bf{D}}}
\newcommand{\BZ}{{\bf{Z}}}
\newcommand{\BSigma}{{\mathbf{\Sigma}}}
\newcommand{\bmu}{{\boldsymbol{\mu}}}
\newcommand{\be}{{\bf{e}}}
\newcommand{\bz}{{\bf{z}}}
\newcommand{\by}{{\bf{y}}}

\newcommand{\pzc}{\textcolor{black}}
\newcommand{\ning}{\textcolor{black}}

\title{Synthetic Acceleration Preconditioners for Parametric Radiative Transfer Equations based on Trajectory-Aware  \ning{Reduced-Order} Models}

\begin{document}

\author{
Ning Tang, Zhichao Peng
}

\ead{ntangad@connect.ust.hk,pengzhic@ust.hk}

\affiliation{organization={Department of Mathematics, The Hong Kong University of Science and Technology},
            addressline={Clear Water Bay, Kowloon},
            city={Hong Kong},
            country={China}
}

\begin{abstract}



The parametric radiative transfer equation (RTE) arises in multi-query applications, such as design optimization, inverse problems, and uncertainty quantification, which require solving the RTE multiple times for various parameters. \ning{Classical synthetic acceleration (SA) preconditioners are typically derived from low-order approximations of the kinetic correction equation, such as the diffusion approximation in diffusion synthetic acceleration (DSA).} 
Despite their widespread success, these methods rely on empirical physical assumptions and do not exploit low-rank structures across parameters of the parametric problem. 

To address this limitation, our previous work developed ROMSAD, a reduced-order model (ROM) enhanced SA preconditioner that incorporates low-rank structures across parameters and the original kinetic description of the correction equation. While ROMSAD improves overall efficiency compared with DSA, its efficiency decreases after the first iteration because the construction of the underlying ROM ignores the preconditioner-dependence of the residual trajectory, leading to a mismatch between the offline and online residual trajectories. 

\ning{To address this offline–online mismatch, we adapt multispace reduced basis preconditioning to parametric RTEs, developing a trajectory-aware ROM framework. Specifically, we sequentially construct iteration-dependent reduced bases to account for the dependence of residual trajectories on the ROM-based preconditioners applied in preceding iterations. We also formulate SA-based kinetic corrections tailored to matrix-free transport sweeps in RTE solvers to achieve efficient online solves.} Numerical experiments demonstrate superior efficiency over DSA and ROMSAD. For a parametric lattice problem, the proposed
preconditioners lead to convergence within only \ning{average $3.4$} transport sweeps online.

\end{abstract}


\maketitle

\section{Introduction}
The radiative transfer equation (RTE) describes how particles propagate through and interact with participating media. It has broad applications, including medical imaging \cite{abdoulaev2003three}, 
astrophysics \cite{rybicki2024radiative}, 
remote sensing \cite{spurr2001linearized}, 
and nuclear engineering \cite{pomraning1973equations}. 
Multi-query applications, such as sensitivity analysis, design optimization, and inverse problems, often require solving the RTE repeatedly for various parameters such as geometric configurations, material properties, and boundary conditions. Consequently, developing efficient iterative solvers becomes crucial for such applications.

Source iteration with synthetic acceleration (SI-SA) is a popular iterative solver for  the RTE \cite{Adams2002FastIM, azmy2010advances}. Each iteration of the SI-SA begins with an SI step. During this SI step, the macroscopic density is frozen to decouple the microscopic particle distribution for various angular directions. The subsequent SA step then accelerates convergence by acting as a preconditioner that corrects this macroscopic density. The ideal correction that ensures convergence of the subsequent iteration is the solution to an ideal correction equation, which is also an RTE. In practice, this ideal correction equation is replaced by a computationally feasible low-order approximation. Diffusion Synthetic Acceleration (DSA) uses its diffusion limit \cite{kopp1963synthetic,reed1971effectiveness,alcouffe1977diffusion}. $S_2$ Synthetic Acceleration ($\mathrm{S}_2$SA) employs an $\mathrm{S}_2$ angular approximation \cite{lorence1989s}. 
The Quasi-Diffusion method leverages the variable Eddington factor \cite{gol1964quasi,anistratov1993nonlinear,olivier2023family}. Furthermore, extensions of SA preconditioners to Krylov methods \cite{warsa2004krylov} are introduced to further improve the efficiency and robustness for problems with strong material interface discontinuities.

Despite the great success of classical SA preconditioners, several limitations remain. First, classical SA preconditioners rely on the effectiveness of the underlying empirical low-order approximation. For example, DSA may become less efficient if the problem is far from its diffusion limit \cite{ren2019fast}. Moreover, they also do not leverage low-rank structures of the solution manifold across parameters of the underlying parametric problem. 

To address these limitations, reduced-order model (ROM) enhanced SA preconditioners, called ROMSAD, have been designed \cite{peng2024romsad,peng2025flexible}. Before delving into details, we briefly review the literature on ROMs for the RTE. Projection-based ROMs typically follow an offline-online decomposition framework. In the offline stage, a low-dimensional basis is built by exploring low-rank structures in the solution data for parametric problems.  In the online stage, great computational savings can be achieved through projection onto the low-dimensional space constructed offline.  Recently, ROMs for the RTE have been actively developed based on proper orthogonal decomposition \cite{buchan2015pod,hughes2020discontinuous,choi2021space,coale2021reduced,tano2021affine,buchan2024reduced,hardy2024proper}, greedy algorithms \cite{tencer2016reduced,peng2022reduced,peng2024micro,matsuda2025reduced}, dynamic mode decomposition \cite{mcclarren2019calculating,smith2023variable}, and proper generalized decomposition \cite{dominesey2022reduced}.

In ROMSAD \cite{peng2024romsad,peng2025flexible}, we first utilize a projection-based ROM as a low-rank approximation to the kinetic correction equation to eliminate low-frequency errors, and then combine it with DSA to handle remaining high-frequency errors. This methodology allows us to exploit the original kinetic description of the ideal correction equation and low-rank structures across parameters. Despite achieving significant online acceleration over DSA for parametric problems, the ROMSAD preconditioner suffers from a reduction in efficiency after the first iteration.

The cause of this efficiency reduction is as follows. The right-hand side of the ideal correction equation for each iteration is determined by the residual for the current iteration. In the offline stage, \cite{peng2024romsad,peng2025flexible} construct the ROM based on residual trajectories generated by a fixed preconditioner, e.g., DSA, while the online residuals are determined dynamically by the ROM itself. This inconsistency prevents the ROM from accurately capturing the iteration-dependent trajectory evolution, leading to progressive efficiency reduction.

To eliminate this mismatch between the offline and online residual trajectories, \ning{we adapt the multispace reduced-basis (MSRB) preconditioning strategy \cite{santo2018multi} to parametric RTEs, developing a trajectory-aware, ROM-enhanced SA framework. During ROM construction, we explicitly account for how previously applied ROM-based preconditioners alter subsequent residual trajectories.} The key components of our method are as follows. 
\begin{enumerate}
    \item In the offline stage, we sequentially construct iteration-dependent ROMs. Each new ROM reflects how the residual trajectory depends on the ROM-based correction from the earlier iterations, thereby removing the offline-online mismatch. \ning{Crucially, we tailor this construction to parametric RTEs so that it integrates seamlessly with matrix-free transport sweeps, enabling efficient online solves.}
    
    \item Similar to \cite{peng2024romsad}, \ning{in constructing} these reduced-order spaces, \ning{we generate} correction snapshots without directly solving the ideal correction equation. When the underlying iterative solver is the flexible generalized minimal residual (FGMRES) method, we exploit a reformulation of the kinetic correction equation \cite{peng2025flexible} to generate the required snapshots.

    \item In the online stage, \ning{ROM-based SA preconditioners are applied} in the first few iterations, followed by DSA in subsequent iterations. This \ning{hybrid} strategy uses ROMs to reduce the dominant error components and DSA to handle the remaining errors.
\end{enumerate}
The key advance of the trajectory-aware framework over ROMSAD lies not merely in using multiple ROMs across iterations, but in its sequential and trajectory-aware ROM construction to mitigate the mismatch between offline and online residual trajectories. \ning{Our method shares the multispace philosophy of MSRB preconditioning \cite{santo2018multi}, but it is specifically tailored to parametric RTEs. Directly applying MSRB to the sweep-based scalar-flux formulation of the RTE with parametric cross sections generally requires additional online transport sweeps to assemble the reduced operator \cite{behne2022minimally}, potentially offsetting the online speedup or even making the method slower than DSA. In contrast, our kinetic correction formulation enables complete offline precomputation of the reduced matrices, eliminating these additional online sweeps.}

In our numerical tests, the proposed trajectory-aware approach achieves significant online acceleration and enhanced robustness over ROMSAD preconditioners \cite{peng2024romsad,peng2025flexible} by paying marginal additional offline costs. 

To better contextualize our approach, we briefly review other acceleration techniques for the RTE based on data-driven ROMs or low-rank approximations. \ning{Several closely related studies \cite{mcclarren2018acceleration,mcclarren2022data,wang2024dynamic} develop acceleration methods based on dynamic mode decomposition (DMD). Specifically, DMD exploits low-rank structure across successive iterations of a single solve by approximating the source-iteration dynamics to identify slowly decaying error modes. Our method is designed to leverage low-rank structures across parameters by constructing iteration-dependent POD spaces from trajectory-aware kinetic correction snapshots. The two approaches therefore have complementary strengths.} A ROM-based nested iteration method is developed in \cite{dahmen2020adaptive}, while randomized singular value decomposition is employed to obtain a low-rank boundary-to-boundary map in a Schwarz solver for the RTE \cite{chen2020random}. Offline–online decomposition has also been applied to accelerate the tailored finite point method \cite{fu2024fast}. The approaches in \cite{dahmen2020adaptive,chen2020random,fu2024fast} do not focus on enhancing SA preconditioners, which is the primary focus of our work. Beyond data-driven ROMs, another line of research considers low-rank solvers based on matrix or tensor decompositions, including preconditioned low-rank iterative solvers \cite{bachmayr2026low,guo2026inexact,haut2026efficient,guo2026highly}, time marching via dynamical low-rank algorithms \cite{einkemmer2021asymptotic,peng2020low,yin2024towards,einkemmer2024asymptotic}, and step-and-truncation methods \cite{sands2024high}. A key advantage of these low-rank solvers is that they are offline-free, while data-driven approaches may offer greater online acceleration and can be combined with them \pzc{as in \cite{shang2026incremental}}.

The rest of the paper is organized as follows. In Sec.~\ref{sec:background}, we introduce the steady-state linear RTE and its corresponding discretization scheme, followed by a brief review of ROMs. Sec.~\ref{sec:source_iteration_ac} reviews SI and ROMSAD, an existing ROM-enhanced SA preconditioner, and identifies their limitations to motivate our trajectory-aware approach. In Sec.~\ref{sec:trajectory_propose}, we develop the trajectory-aware framework and combine it with both SI and Krylov methods. Performance of our method is demonstrated through numerical experiments in Sec.~\ref{sec:numerical}, while conclusions are drawn in Sec.~\ref{sec:conclusion}.


\section{Governing Equation and its discretization scheme\label{sec:background}}
The radiative transfer equation (RTE) serves as the governing equation for characterizing the propagation and interaction of particle systems within participating media. In this paper, we focus on the single energy group steady-state linear RTE with isotropic scattering and isotropic inflow boundary conditions:
\begin{subequations}
\label{eq:rte}
    \begin{align}
    &\Bupsilon\cdot\nabla_{\mathbf{r}}f(\mathbf{r},\Bupsilon) + \sigma_{t}(\mathbf{r})f(\mathbf{r},\Bupsilon) = \sigma_{s}(\mathbf{r})\rho(\mathbf{r}) + Q(\mathbf{r}),\;\; \rho(\mathbf{r}) = \frac{1}{4\pi}\int_{\mathbb{S}^2} f(\mathbf{r},\Bupsilon) d\Bupsilon, \;\; \mathbf{r}\in \mathcal{D}, \;\; \Bupsilon\in\mathbb{S}^2,\label{eq:governing}
    \\
    &f(\mathbf{r}, \Bupsilon) = g(\mathbf{r}), \;\;\mathbf{r}\in\partial\mathcal{D},\;\; \Bupsilon\cdot \mathbf{n}(\mathbf{r})<0.
    \label{eq:boundary}
    \end{align}
\end{subequations}
Here, $f(\mathbf{r}, \Bupsilon)$ is the distribution function (also known as the angular flux) at position $\mathbf{r}\in\mathcal{D}\subset \mathbb{R}^n$ and in the angular direction $\Bupsilon\in\mathbb{S}^2$, where $\mathbb{S}^2$ is the unit sphere. $\rho(\mathbf{r})$ is the macroscopic density (also known as the scalar flux). $Q(\mathbf{r})$ is an isotropic source, and $\sigma_{t}(\mathbf{r}) = \sigma_{a}(\mathbf{r}) + \sigma_{s}(\mathbf{r})$ is the total cross section combining the absorption cross section $\sigma_a(\br)\geq 0$ and the scattering cross section $\sigma_s(\br)\geq 0$.  This system is closed by the inflow boundary condition \eqref{eq:boundary}, where $\mathbf{n}(\mathbf{r})$ is the outward normal vector field at boundary point $\mathbf{r}\in\partial\mathcal{D}$.

We note that the RTE has multiscale behavior. In the free-streaming limit with $\sigma_t=0$, it is a pure transport equation. When the scattering cross section $\sigma_{s}\rightarrow \infty$, it asymptotically converges to its diffusion limit \cite{bardos1984diffusion}.

In this section, we start with the numerical discretization for the RTE using the discrete ordinates $(S_{N})$ method and the upwind discontinuous Galerkin method, followed by a brief review of ROMs.

\subsection{Spatial-angular discretization\label{sec:upwind dg}}
\textbf{Discrete ordinates ($S_{N}$):} The $S_{N}$ method approximates the macroscopic density by solving the equation only at a set of  quadrature points for $\mathbb{S}^2$, namely $\{\Bupsilon_{j}\}_{j=1}^{N_{\Bupsilon}}$,  with corresponding quadrature weights $\{\omega_{j}\}_{j=1}^{N_{\Bupsilon}}$. The integral term $\rho(\mathbf{r})$ can then be approximated by
\begin{equation}
    \rho(\mathbf{r})  = \frac{1}{4\pi}\int_{\mathbb{S}^2} f(\mathbf{r},\Bupsilon) d\Bupsilon \approx \sum_{j=1}^{N_{\Bupsilon}} \omega_{j} f(\mathbf{r},\Bupsilon_{j}), \;\; \sum_{j=1}^{N_{\Bupsilon}}\omega_{j} = 1.
\end{equation}
Thus, the RTE \eqref{eq:rte} can be reduced to a set of partial differential equations:
\begin{subequations}
\label{eq:rte discrete}
    \begin{align}
    &\Bupsilon_{j}\cdot\nabla_{\mathbf{r}}f(\mathbf{r},\Bupsilon_{j}) + \sigma_{t}(\mathbf{r})f(\mathbf{r},\Bupsilon_{j}) = \sigma_{s}(\mathbf{r})\rho(\mathbf{r}) + Q(\mathbf{r}),\;\; \rho(\mathbf{r}) = \sum_{j=1}^{N_{\Bupsilon}} \omega_{j} f(\mathbf{r},\Bupsilon_{j}), \label{eq:governing_discrete}
    \\
    &f(\mathbf{r}, \Bupsilon_{j}) = g(\mathbf{r}), \;\;\mathbf{r}\in\partial\mathcal{D},\;\; \Bupsilon_{j}\cdot \mathbf{n}(\mathbf{r})<0.
    \label{eq:boundary_discrete}
    \end{align}
\end{subequations}
Here, we employ different quadrature rules for different geometric configurations. For 1D slab geometry, we utilize the Gauss-Legendre quadrature. In higher dimensions, we use Chebyshev-Legendre (CL) quadratures. The CL quadrature rule, CL($N_{\alpha}, N_{\Bupsilon_z}$), is constructed as the tensor product of two normalized quadrature schemes: the $N_{\alpha}$-point Chebyshev rule for the unit circle 
\begin{equation}
\left\{(\alpha_j,\omega_{\alpha,j}): \;\alpha_j = \frac{(2j-1)\pi}{N_\alpha}, \;\;\omega^\alpha_j=\frac{1}{N_\alpha},\; j=1,\dots,N_\alpha\right\}
\end{equation}
and the $N_{\Bupsilon_z}$-point Gauss-Legendre quadrature rule $\{(\Bupsilon_{z,j},\omega_{z,j})\}_{j=1}^{N_{\Bupsilon_z}}$ for the $z$-component of the angular direction  $\Bupsilon_z\in[-1,1]$, where $\Bupsilon_{z,j}$ denotes the \(j\)-th root of the Legendre polynomial of degree $N_{\Bupsilon_{z}}$.
Then, the quadrature points and related weights of the CL($N_{\alpha}$,$N_{\Bupsilon_z}$) quadrature rule are given by
\begin{equation}
\Bupsilon_{j}=\left(\cos(\alpha_{j_1})\sqrt{1-\Bupsilon_{z,j_2}^2},\;\;\sin(\alpha_{j_1})\sqrt{1-\Bupsilon_{z,j_2}^2},\;\;\Bupsilon_{z,j_2}\right)\quad,\;\;\omega_j=\omega_{\alpha,j_1}\omega_{z,j_2},
\end{equation}
where  $1\leq j_1\leq N_\alpha$, $1\leq j_2\leq N_{\Bupsilon_z}$, $1\leq j=(j_2-1) N_{\alpha}+j_1\leq N_{\Bupsilon}$, and in this case $N_{\Bupsilon}=N_{\alpha}N_{\Bupsilon_z}$ for \eqref{eq:governing_discrete}.

\textbf{Upwind discontinuous Galerkin method (DG):} 
For spatial discretization, we apply a high-order upwind discontinuous Galerkin (DG) method, which is proven to be \ning{asymptotic-preserving (AP)} \cite{adams2001discontinuous,guermond2010asymptotic}. We consider a 2D rectangular computational domain and partition it with a family of rectangular elements $\mathcal{I}_{h} = \{\mathcal{I}_{i}\}_{i=1}^{N_{\mathbf{r}}}$. The DG finite element space is defined as
\begin{equation}
    \mathrm{U}_{h}^{K} = \left\{ u\in \mathrm{L}^2(\mathcal{D}): u|_{\mathcal{I}_i}\in \mathbb{Q}^{K}(\mathcal{I}_{i}),\;\;\forall\;\mathcal{I}_{i}\in\mathcal{I}_{h}\right\},
    \label{eq:dg_space}
\end{equation}
where $\mathbb{Q}^{K}$
is the set of polynomials defined on $\mathcal{I}_i$ whose degree in each direction is at most
$K$. Let $\partial \mathcal{I}_{h}$ denote the set of cell edges, and
\begin{equation}
    \partial\mathcal{I}_{h,j}^{bc} = \left\{ \mathcal{E}\in \partial\mathcal{D} \cap\partial\mathcal{I}_h: \Bupsilon_{j}\cdot\mathbf{n}(\mathcal{E})<0 \right\}
\end{equation}
is the subset of the edges where inflow boundary conditions for the angle  $\Bupsilon_{j}$ are imposed. Then, the DG spatial discretization with upwind flux for the $S_N$
system \eqref{eq:governing_discrete} is to seek $f_{h}(\mathbf{r},\Bupsilon_{j})\in\mathrm{U}_{h}^{K}, j =1,2,...,N_{\Bupsilon}$ such that  $\forall\phi_{h}(\mathbf{r})\in\mathrm{U}_{h}^{K}$,
\begin{align}
\displaystyle
-\sum_{i=1}^{N_{\mathbf{r}}}&\int_{\mathcal{I}_i} \Big(\Bupsilon_j\cdot\nabla\phi_h(\mathbf{r})\Big) f_h(\mathbf{r},\Bupsilon_j) d\mathbf{r}+\sum_{\mathcal{E}\in\partial\mathcal{I}_h\setminus\partial\mathcal{I}_{h,j}^{\textrm{bc}}} \int_{\mathcal{E}} \widehat{\BW}(\Bupsilon_j,f_h, \bn(\mathbf{r}))\phi_h(\mathbf{r}) d\mathbf{r}+\sum_{i=1}^{N_{\mathbf{r}}}\int_{\mathcal{I}_i}\sigma_t(\mathbf{r}) f_h(\mathbf{r},\Bupsilon_j)\phi_h(\mathbf{r}) d\mathbf{r}
\notag
\\
= &\sum_{i=1}^{N_{\mathbf{r}}}\int_{\mathcal{I}_i}\sigma_s(\mathbf{r}) \rho_h(\mathbf{r})\phi_h(\mathbf{r}) d\mathbf{r}
 + \sum_{i=1}^{N_{\mathbf{r}}}\int_{\mathcal{I}_i}Q(\mathbf{r})\phi_h(\mathbf{r}) d\mathbf{r}
 -\sum_{\mathcal{E}\in\partial\mathcal{I}_{h,j}^{\textrm{bc}}}\int_{\mathcal{E}} g(\mathbf{r})\phi_h(\mathbf{r}) \Bupsilon_j\cdot \bn(\mathbf{r}) d\mathbf{r},
 \label{eq:DG}
\end{align}
where the upwind flux $\widehat{\BW}(\Bupsilon_j,f_h, \bn(\mathbf{r}))$ along the interface $\mathcal{E}$ is defined as
\begin{align}
\widehat{\BW}(\Bupsilon_j, f_h,\bn(\mathbf{r}))\Big|_{\mathcal{E}} = \frac{\Bupsilon_j\cdot\bn(\mathbf{r}) + |\Bupsilon_j\cdot \bn(\mathbf{r})|}{2}f_h^-(\mathbf{r},\Bupsilon_j)+\frac{\Bupsilon_j\cdot\bn(\mathbf{r}) -|\Bupsilon_j\cdot \bn(\mathbf{r})|}{2}f_h^+(\mathbf{r},\Bupsilon_j).
\label{eq:upwind}
\end{align}
Here, the superscripts ``$-$'' and ``$+$'' denote the values taken from the two neighboring elements $\mathcal I^-$ and $\mathcal I^+$, respectively, and $\bn(\mathbf r)$ denotes the outward unit normal vector of $\mathcal I^-$. 

\subsection{Matrix-Vector formulation\label{sec:matrix}}
Let \(\{\phi_i\}_{i=1}^{N_{\mathrm{DOF}}}\) be an orthonormal basis for the DG space. Using the weak formulation \eqref{eq:DG} and the upwind flux \eqref{eq:upwind}, we obtain the following fully discrete scheme for the RTE \eqref{eq:governing}:
\begin{align}
    (\BD_j+\BSigma_t) \bff_j =\BSigma_s\brho+\bQ+\bg_j^{(\textrm{bc})}=\BSigma_s\brho+\widetilde{\bQ}_j,\;\;
\brho=\sum_{j=1}^{N_{\Bupsilon}}\omega_j\bff_j,\;\;j=1,\dots,N_{\Bupsilon},
\label{eq:dg_matrix_vec}
\end{align}
where $\bff_{j}\in\mathbb{R}^{N_{\textrm{DOF}}}$, $\brho\in\mathbb{R}^{N_{\textrm{DOF}}}$ are the vectors collecting degrees of freedom (DOFs) for $f_{h}(\cdot,\Bupsilon_j)$ and $\rho_h$, respectively. The discrete advection operator and cross-section matrices $\BD_j,\BSigma_t,\BSigma_s\in\mathbb{R}^{N_{\textrm{DOF}}\times N_{\textrm{DOF}}}$, and the source and boundary flux $\bQ,\bg_j^{(\textrm{bc})}\in\mathbb{R}^{N_\textrm{DOF}}$ are defined as:
\begin{subequations}
\begin{align}
&(\BD_{j})_{mn} = -\sum_{i=1}^{N_{\mathbf{r}}}\int_{\mathcal{I}_i} (\Bupsilon_j\cdot\nabla\phi_m(\mathbf{r})) \phi_n(\mathbf{r}) d\mathbf{r}+ \sum_{\mathcal{E}\in\partial\mathcal{I}_h\setminus\partial\mathcal{I}_{h,j}^{\textrm{bc}}}\int_{\mathcal{E}} \widehat{\BW}\left(\Bupsilon_j, \phi_n,\bn(\mathbf{r})\right)\phi_m(\mathbf{r}) d\mathbf{r},
\\
&(\BSigma_s)_{mn} = \sum_{i=1}^{N_{\mathbf{r}}}\int_{\mathcal{I}_i} \sigma_s(\mathbf{r})\phi_m(\mathbf{r}) \phi_n(\mathbf{r}) d\mathbf{r},\quad (\BSigma_t)_{mn} = \sum_{i=1}^{N_{\mathbf{r}}}\int_{\mathcal{I}_i} \sigma_t(\mathbf{r})\phi_m(\mathbf{r}) \phi_n(\mathbf{r}) d\mathbf{r},
\label{eq:element_sigma}
\\
&(\bQ)_m= \sum_{i=1}^{N_{\mathbf{r}}}\int_{\mathcal{I}_i} Q(\mathbf{r})\ning{\phi}_m(\mathbf{r}) d\mathbf{r},\quad
(\bg_j^{(\textrm{bc})})_m = -\sum_{\mathcal{E}\in\partial\mathcal{I}_{h,j}^{\textrm{bc}}}\int_{\mathcal{E}} g(\mathbf{r})\ning{\phi_m}(\mathbf{r}) \Bupsilon_j\cdot \bn(\mathbf{r}) d\mathbf{r}.
\end{align}
\end{subequations}
Equation \eqref{eq:dg_matrix_vec} can be written in block-matrix form as:
\begin{equation}
\label{eq:one_equation}
\BA \mathbf{f} =
\left(\begin{matrix}
\BD_1+\BSigma_t - \omega_1\BSigma_s & -\omega_2 \BSigma_s & \dots & -\omega_{N_{\Bupsilon}}\BSigma_s\\
-\omega_1\BSigma_s & \BD_2+\BSigma_t-\omega_2\BSigma_s & \dots & -\omega_{N_{\Bupsilon}}\BSigma_s \\
\vdots & \vdots & \vdots & \vdots \\
-\omega_1\BSigma_s & -\omega_2\BSigma_s & \dots & \BD_{N_{\Bupsilon}}+\BSigma_t-\omega_{N_{\Bupsilon}}\BSigma_s
\end{matrix}\right)
\left(
\begin{matrix}
\bff_1\\
\bff_2\\
\vdots\\
\bff_{N_{\Bupsilon}}
\end{matrix}
\right)
=
\left(\begin{matrix}
\widetilde{\bQ}_1\\
\widetilde{\bQ}_2\\
\vdots\\
\widetilde{\bQ}_{N_{\Bupsilon}}
\end{matrix}
\right)=\bb.
\end{equation}

\subsection{Reduced-order model\label{sec:rom}}
Multi-query applications, e.g., sensitivity analysis, inverse problems, and uncertainty quantification, require solving the RTE for various parameters such as material properties, configurations, and boundary conditions. A reduced-order model (ROM) is a framework to achieve fast online computations for such parametric cases.

Most ROM approaches follow an offline-online decomposition framework. During the offline stage, a low-dimensional subspace is constructed to extract low-rank features of the underlying parametric problem from high-fidelity solution data. The online stage enables rapid solution computation by projecting the full-order equations onto this low-dimensional subspace constructed offline or performing interpolation in it.

Let $\mathscr{P}$ denote a set of parameters and $\mathscr{P}_{\mathrm{train}} = \{\bmu_i\}_{i=1}^{N_{\mathrm{train}}}$ denote the corresponding training set. For the RTE, these parameters could represent material properties such as scattering and absorption cross sections. With this parameter dependence, the governing equation \eqref{eq:one_equation} can be denoted as
\begin{equation}
    \mathbf{A}_\bmu \bff_\bmu = \mathbf{b_\bmu}.
    \label{eq:parameter_equation}
\end{equation}
Here, we assume that both $\mathbf{A}_\bmu$ and $\mathbf{b_\bmu}$ admit affine decompositions: 
\begin{equation}
    \mathbf{A}_\bmu = \sum_{q=1}^{k_a} \theta^{\mathbf{A}}_{q}(\bmu)\mathbf{A}_q, \;\; \mathbf{b}_\bmu = \sum_{p=1}^{k_b} \theta^{\mathbf{b}}_{p}(\bmu)\mathbf{b}_p,
\end{equation}
where $\theta^{\mathbf{A}}_{q}, \theta^{\mathbf{b}}_{p}: \mathscr{P}\rightarrow \mathbb{R}$ are  scalar functions, and  $\mathbf{A}_q\in \mathbb{R}^{N_{h}\times N_{h}}$, $\mathbf{b}_p \in\mathbb{R}^{N_{h}}$ are independent of parameter $\bmu$. 

\textbf{Offline construction:} For each parameter vector $\bmu_i\in\mathscr{P}_{\mathrm{train}}$, we denote its corresponding high-fidelity solution as $\bff_{\bmu_i}\in\mathbb{R}^{N_{h}}$, where $N_h = N_{\Bupsilon} N_{\mathrm{DOF}}$. We assemble the snapshot matrix  $\mathbf{F}\in \mathbb{R}^{N_{h}\times N_{\mathrm{train}}}$ as follows:
\begin{equation}
    \mathbf{F} = \left[\bff_{\bmu_1}, \bff_{\bmu_2},\cdots,\bff_{\bmu_{N_{\mathrm{train}}}}\right].
    \label{eq:snapshot}
\end{equation}
To build a ROM, we compute the singular value decomposition (SVD) of the snapshot matrix:
\begin{equation}
\mathbf{F} = \mathbf{U\Sigma V}^{T},
\end{equation}
where $\mathbf{U}\in\mathbb{R}^{N_{h}\times N_{h}}$ contains the proper orthogonal modes, $\mathbf{V}\in\mathbb{R}^{N_{\textrm{train}}\times N_{\textrm{train}}}$, and the $k$-th diagonal element of $\mathbf{\Sigma}\in\mathbb{R}^{N_h\times N_{\textrm{train}}}$ is the $k$-th singular value of $\BF$, namely $s_k$. The dominant low-rank structure is captured by the reduced basis, which consists of the first $r$ columns of $\mathbf{U}$. Here, $r$ is the smallest integer satisfying the following criterion for the threshold $\epsilon_{\textrm{svd}}$: 
\begin{equation}
\mathbf{U}_r = \mathbf{U}( \;: \;,1:r),\quad\textrm{where}\quad\frac{\sum_{i=1}^{r} \ning{s_i^2}}{\sum_{i=1}^{N_{\mathbf{sv}}} \ning{s_i^2}} \geq 1-\epsilon_{\mathrm{svd}}^{\ning{2}},\quad N_{\mathbf{sv}}=\min(N_h,N_{\textrm{train}}).
\end{equation}

\textbf{Online projection:} For any new parameter configuration $\bmu\notin\mathscr{P}_{\mathrm{train}}$, we approximate $\bff_\bmu$ by a low-rank solution in the form $\bff_{\bmu}\approx\bff_{\bmu,r}=\mathbf{U}_{r}\bc_{\bmu,r}$, where $\bc_{\bmu,r}\in\mathbb{R}^r$. Here,  $\bc_{\bmu,r}$ is determined by solving an $r$-dimensional system $\mathbf{A}_{\bmu,r}\mathbf{c}_{\bmu,r} = \mathbf{b}_{\bmu,r}$ derived through the projection of equation
\eqref{eq:parameter_equation}, where
\begin{subequations}
\label{eq:reduced_equation}
    \begin{align}
    & \mathbf{A}_{\bmu,r} = \mathbf{U}_{r}^T\mathbf{A}_\bmu\mathbf{U}_r = \sum_{q=1}^{k_a}\theta^{\mathbf{A}}_{q}(\bmu)\mathbf{U}_{r}^T\mathbf{A}_q\mathbf{U}_{r},\label{eq:precomputed_term_1}\\
    &\mathbf{b}_{\bmu,r} = \mathbf{U}_{r}^T\mathbf{b}_{\bmu} = \sum_{p=1}^{k_b} \theta^{\mathbf{b}}_{p}(\bmu)\mathbf{U}_{r}^T\mathbf{b}_p.
    \label{eq:precomputed_term_2}
    \end{align}
\end{subequations}

\textbf{Offline precomputation:} $\BA_{\bmu,r}$ and $\bb_{\bmu,r}$ can be efficiently assembled online using the affine decomposition \eqref{eq:reduced_equation} and offline-precomputed matrices $\mathbf{U}_{r}^T\mathbf{A}_q\mathbf{U}_{r} \in\mathbb{R}^{r\times r}$ and vectors $\BU_r^T\bb_{p}$. Moreover, the numerical integration to compute density can also be accelerated through offline precomputations. Let $\mathbf{U}_{r,j}$ denote the block of rows of $\mathbf{U}_r$ aligned with $\bff_{j}$ (the $( (j-1)N_{\mathrm{DOF}}+1) \mathrm{-th}$ row to the $(jN_{\mathrm{DOF}})\mathrm{-th}$ row). Then, the density can be efficiently approximated online by:
\begin{equation}
    \brho_{\bmu} \approx \brho_{\bmu,r}=\sum_{j=1}^{N_{\Bupsilon}}
\omega_j\mathbf{U}_{r,j}\bc_{\bmu,r} = \mathbf{U}_r^{\rho}\bc_{\bmu,r}\;,
\end{equation}
where $\mathbf{U}_r^{\rho} = \sum_{j=1}^{N_{\Bupsilon}}
\omega_j\mathbf{U}_{r,j}$ is precomputed offline. 

\begin{rem}
When the parameterized operator $\mathbf{A}_\bmu$ exhibits non-affine parametric dependence, hyper-reduction techniques, e.g., Empirical Interpolation Method (EIM) \cite{barrault2004empirical} and Discrete EIM (DEIM) \cite{chaturantabut2010nonlinear}, can be employed to maintain the online efficiency.
\end{rem}

\section{Source iteration and synthetic acceleration \label{sec:source_iteration_ac}}
In this section, we provide an overview of classical synthetic acceleration (SA) preconditioners for source iteration (SI) and Krylov methods, followed by ROMSAD, our previous ROM-enhanced SA preconditioner \cite{peng2024romsad,peng2025flexible}. We motivate our new trajectory-aware approach by pointing out the limitations of these methods. 

\subsection{Source iteration and synthetic acceleration}
We first outline SI-SA, followed by a detailed algorithm in Alg. \ref{alg:sisa}.

\textbf{Source iteration (SI):} 
SI decouples the angular directions by freezing the density $\brho$. In the $l$-th iteration, the following decoupled subsystems are solved:
\begin{equation}
    (\BD_j+\BSigma_t) \bff_j^{(l)} =\BSigma_s\brho^{(l-1)}+\widetilde{\bQ}_j,\;\; j=1,\dots,N_{\Bupsilon}, \;\;\forall\;l\geq 1, \label{eq:source_iteration}
\end{equation}
where $\brho^{(0)}$ is obtained from the initial guess. With upwind discretization, the operator $\BD_j+\BSigma_t$ is block lower triangular after reordering elements along the upwind direction of $\Bupsilon_j$. This structure enables efficient matrix-free solves using transport sweeps \cite{adams2001discontinuous,azmy2010advances}. After obtaining  $\bff_j^{(l)}$ for all angular directions, the macroscopic density is updated as: $\brho^{(l,*)}=\sum_{j=1}^{N_{\Bupsilon}}\omega_j\bff_j^{(l)}$.

\textbf{Synthetic acceleration (SA):} Without preconditioning, SI may suffer from slow convergence or even false convergence \cite{Adams2002FastIM}. SA serves as  a preconditioning step to address this issue by introducing a density correction $\delta\brho$ after each SI step: 
\begin{equation}
\brho^{(l)} = \brho^{(l,*)} + \delta\brho^{(l)},
\;\;\textrm{where}\;\;\brho^{(l,*)}=\sum_{j=1}^{N_{\Bupsilon}}\omega_j\bff_j^{(l)}.
\end{equation}
Taking $\bff_j$ as the exact solution of \eqref{eq:dg_matrix_vec}, we then define the ideal density correction as follows:
\begin{equation}
    \delta \bff_{j}^{(l)} = \bff_j - \bff_j^{(l)}, \;\; \delta\brho^{(l)} = \sum_{j=1}^{N_{\Bupsilon}}\omega_j\delta\bff_j^{(l)}.
    \label{eq:ideal_density}
\end{equation}
By subtracting \eqref{eq:source_iteration} from \eqref{eq:dg_matrix_vec}, we obtain the discrete ideal correction equation:
\begin{equation}
     (\BD_j+\BSigma_t) \delta\bff_j^{(l)} =\BSigma_s\delta\brho^{(l)} + \BSigma_s(\brho^{(l,*)}-\brho^{(l-1)}),\;\; j=1,\dots,N_{\Bupsilon}. \;\; \label{eq:ideal_density_equation}
\end{equation}
Following the notation in \eqref{eq:ideal_density}, \eqref{eq:ideal_density_equation} can be written as $\mathbf{A}\delta\mathbf{f}^{(l)} = \delta\mathbf{b}^{(l)}$, 
where $\mathbf{A}\in\mathbb{R}^{N_h\times N_h}$ is defined in \eqref{eq:one_equation}, and
\begin{subequations}
    \begin{align}
    & \delta\mathbf{f}^{(l)} = \left((\delta \bff_{1}^{(l)})^{T},(\delta \bff_{2}^{(l)})^{T},\;\cdots\;, (\delta \bff^{(l)}_{N_{\Bupsilon}} )^{T} \right)^{T} \in \mathbb{R}^{N_{h}},
    \\
    & \delta\mathbf{b}^{(l)} = \left((\BSigma_s(\brho^{(l,*)}-\brho^{(l-1)}))^T ,\; \cdots\;,(\BSigma_s(\brho^{(l,*)}-\brho^{(l-1)}))^T\right)^{T}  \in\mathbb{R}^{N_{h}}.
    \end{align}
\end{subequations}
The correction equation \eqref{eq:ideal_density_equation} 
 can be seen as a discretization of the following RTE with an isotropic source term $\sigma_s(\mathbf{r})(\rho^{(l,*)}(\mathbf{r}) - \rho^{(l-1)}(\mathbf{r}))$ and zero inflow boundary conditions:
\begin{subequations}
\label{eq:correction}
    \begin{align}
    &\Bupsilon\cdot\nabla_{\mathbf{r}}\delta f^{(l)}(\mathbf{r},\Bupsilon) + \sigma_{t}(\mathbf{r})\delta f^{(l)}(\mathbf{r},\Bupsilon) = \sigma_{s}(\mathbf{r})\delta\rho^{(l)}(\mathbf{r}) + \sigma_{s}(\mathbf{r})(\rho^{(l,*)}(\mathbf{r}) - \rho^{(l-1)}(\mathbf{r})),\;\;\mathbf{r}\in \mathcal{D}, \\
    &\delta\rho^{(l)}(\mathbf{r}) = \frac{1}{4\pi}\int_{\mathbb{S}^2} \delta f^{(l)}(\mathbf{r},\Bupsilon) d\Bupsilon, \;\; \mathbf{r}\in \mathcal{D}, \;\; \Bupsilon\in\mathbb{S}^2,
    \\
    &\delta f^{(l)}(\mathbf{r}, \Bupsilon) = 0, \;\;\mathbf{r}\in\partial\mathcal{D},\;\; \Bupsilon\cdot \mathbf{n}(\mathbf{r})<0.
    \end{align}
\end{subequations}

Theoretically, an exact solution of the correction equation \eqref{eq:correction} ensures convergence within at most two SI iterations. However, solving it is computationally equivalent to solving the original RTE, making direct implementation impractical. In practice, the correction equation \eqref{eq:correction} is replaced by a computationally feasible approximation. For example, Diffusion Synthetic Acceleration (DSA) \cite{alcouffe1977diffusion, larsen1982unconditionally} 
solves its diffusion limit: 
\begin{subequations}
\label{eq:diffusion-limit}
\begin{align}
    -\nabla_\mathbf{r}\cdot\left(\frac{1}{\ning{\sigma_t}}\mathbf{D}_{\Bupsilon}(\nabla_\mathbf{r}\delta\rho^{(l)}(\mathbf{r}))\right) + \sigma_a(\mathbf{r})\delta \rho^{(l)}(\mathbf{r}) & = \sigma_{s}(\mathbf{r})(\rho^{(l,*)}(\mathbf{r}) - \rho^{(l-1)}(\mathbf{r})), \quad\ning{\br\in\mathcal{D}},\\
    \ning{\delta\rho^{(l)}
+\frac{2}{3\sigma_t}
\mathbf n\cdot\nabla_{\mathbf r}\delta\rho^{(l)}}
&=0,\quad\br\in\partial\mathcal{D},
\end{align}
\end{subequations}
where $\bn$ is the outward normal direction and $\mathbf{D}_{\Bupsilon} = \textrm{diag}\left(\frac{1}{4\pi}\int_{\mathbb{S}^2} \Bupsilon_x^{2}d\Bupsilon, \;\frac{1}{4\pi}\int_{\mathbb{S}^2} \Bupsilon_y^{2}d\Bupsilon, \;\frac{1}{4\pi}\int_{\mathbb{S}^2} \Bupsilon_z^{2}d\Bupsilon\right).$
Other alternative approximations include the variable Eddington factor in the Quasi-Diffusion method (QD) \cite{gol1964quasi} and low-order $S_N$ angular approximation in $S_2$-SA \cite{lorence1989s} and Transport SA \cite{ramone1997transport}.

\begin{algorithm}[!ht]
\caption{Source iteration and synthetic acceleration (SI-SA). \label{alg:sisa} }
\begin{algorithmic}[1]
\STATE{Given an initial guess $\brho^{(0)}$, tolerance $\epsilon_{\mathrm{SISA}}$, and the maximum number of iterations $N_{\mathrm{iter}}$.}
\FOR{$l = 1:N_{\mathrm{iter}}$}
    \STATE{\textbf{Source iteration:}}
    \STATE{Solve $(\BD_{j}+\BSigma_{t})\bff_{j}^{(l)} =\BSigma_{s}\brho^{(l-1)} + \widetilde{\bQ}_{j}$ with transport sweeps, $j = 1,2,\cdots,N_{\Bupsilon}$.} 
    \STATE{Update the density flux: $\brho^{(l,*)}=\sum_{j=1}^{N_{\Bupsilon}}\omega_j\bff_j^{(l)}.$}
    \IF{the stopping criterion is satisfied}
        \STATE{return solutions: $\brho^{(l,*)}$ and $\bff^{(l)}$}.
    \ELSE
    \STATE{\textbf{Correction:}}
    \STATE{Solve a correction equation to obtain the correction $\delta \brho^{(l)}$}.
    \STATE{Update the density flux: $\brho^{(l)} = \brho^{(l,*)} + \delta \brho^{(l)}$}.
    \ENDIF
\ENDFOR
\end{algorithmic}
\end{algorithm}

\subsection{Extension to Krylov methods \label{sec:sa-krylov}}
As pointed out in \cite{warsa2004krylov}, SI may become ineffective when material properties exhibit strong discontinuities. To address this issue, the SA preconditioner is extended to the Krylov subspace method in \cite{warsa2004krylov} by interpreting SI-SA as a left-preconditioned Richardson iteration. Here, we outline the key idea.

To establish the equivalence between SI-SA and Richardson iteration, we first rewrite \eqref{eq:dg_matrix_vec} as:
\begin{equation}
    \bff_j =(\BD_j+\BSigma_t)^{-1}\BSigma_s\brho+(\BD_j+\BSigma_t) ^{-1}\widetilde{\bQ}_j,\;\; j=1,\dots,N_{\Bupsilon}.
    \label{eq:rewrite_f}
\end{equation}
By applying angular quadrature to \eqref{eq:rewrite_f} and rearranging terms, we obtain:
\begin{equation}
    \left(\mathbf{I} - \sum_{j=1}^{N_{\Bupsilon}} (\omega_j(\BD_j+\BSigma_t)^{-1})\BSigma_s\right)\brho =\sum_{j=1}^{N_{\Bupsilon}} (\omega_j(\BD_j+\BSigma_t) ^{-1})\widetilde{\bQ}_j.
    \label{eq:rewrite_rho}
\end{equation}
Defining $\mathbf{K} = \sum_{j=1}^{N_{\Bupsilon}} (\omega_j(\BD_j+\BSigma_t)^{-1})$ and $\ning{\widetilde{\mathbf{b}} = \sum_{j=1}^{N_{\Bupsilon}} (\omega_j(\BD_j+\BSigma_t) ^{-1})\widetilde{\bQ}_j}$, we rewrite \eqref{eq:rewrite_rho} as
\begin{equation}
    (\mathbf{I} - \mathbf{K}\BSigma_s) \brho = \widetilde{\mathbf{b}}.
    \label{eq:matrix_vec_rewrite_rho}
\end{equation}
Thus, SI without preconditioning can be written as a Richardson iteration for \eqref{eq:matrix_vec_rewrite_rho}:
\begin{equation}
\begin{split}
    \brho^{(l,*)} &= \mathbf{K}\BSigma_s\brho^{(l-1)} + \widetilde{\mathbf{b}}=\brho^{(l-1)} + (\widetilde{\mathbf{b}} - (\mathbf{I} - \mathbf{K}\BSigma_s)\brho^{(l-1)})\\
    & = \brho^{(l-1)} + \boldsymbol{r}^{(l-1)},
    \label{eq:rho_residual}
\end{split}
\end{equation}
where $\boldsymbol{r}^{(l-1)}=\brho^{(l,*)}-\brho^{(l-1)}$ is the residual of equation \eqref{eq:matrix_vec_rewrite_rho} at the $(l-1)$-th step.

In practice,  a computationally tractable low-order operator $\mathbf{C}$ is applied in SA to approximate the kinetic correction operator associated with \eqref{eq:ideal_density_equation}. Following the steps in \ref{sec:left_preconditioner}, the SA step can then be interpreted as introducing a left preconditioner $\BM^{-1}=\BI+\BC^{-1}\BSigma_s$ to \eqref{eq:matrix_vec_rewrite_rho}, i.e.,
\begin{equation}
    (\mathbf{I}+\mathbf{C}^{-1}\BSigma_{s})(\mathbf{I} - \mathbf{K}\BSigma_s) \brho = (\mathbf{I}+\mathbf{C}^{-1}\BSigma_{s})\widetilde{\mathbf{b}}.
    \label{eq:left_preconditioner}
\end{equation}
This viewpoint enables direct extension of SA preconditioners from SI to Krylov methods.


\subsection{Reduced-order model enhanced SI-SA\label{sec:rom_enhanced_sisa}}
The performance of classical SA methods may deteriorate when the accuracy of their underlying low-order approximations is limited \cite{ren2019fast}. Moreover, they do not exploit the underlying low-rank structures of the solution manifold across parameters. Consequently, using a ROM to leverage such low-rank structures to enhance SI-SA naturally emerges for solving parametric problems. With this motivation, a ROM-enhanced SA preconditioner, namely ROMSAD, was proposed for SI \cite{peng2024romsad} and extended to Krylov methods \cite{peng2025flexible}. Its key ideas are as follows.

\begin{enumerate}
    \item ROMSAD combines ROM and DSA to handle different error components. In the early stage of SI, the error is dominated by low-frequency components, which can be effectively captured by ROM. As iteration progresses, high-frequency oscillatory components become dominant. DSA can be seen as a $p$-multigrid applying $P_1$ approximation in the angular space. Thus, it can effectively handle such high-frequency errors. 
    Based on these observations, ROMSAD uses a ROM-based correction during initial iterations and then switches to DSA. See Alg. \ref{alg:romsad} for details.

    \item Solutions of the ideal kinetic correction equation \eqref{eq:ideal_density_equation} are needed to construct the ROM. Instead of solving \eqref{eq:ideal_density_equation} directly, we exploit its definition $\delta \bff^{(l)} = \bff - \bff^{(l)}$. Specifically, when solving the parametric RTE \eqref{eq:source_iteration} with SI-DSA, we store the converged solution $\bff_\bmu$ together with intermediate solutions $\bff_\bmu^{(l)}$ for the first few iterations, $1 \leq l \leq \mathfrak{w}$. The desired correction snapshots are then formed directly from its definition, $\delta \bff^{(l)}_{\bmu} = \bff_{\bmu} - \bff_{\bmu}^{(l)}$. 
\end{enumerate}

\begin{algorithm}[!ht]
\caption{Selection of ROMSA or DSA after the $l$-th source iteration in ROMSAD.\label{alg:romsad}}
\begin{algorithmic}[1]
\STATE \textbf{Input:} Current iteration index $l$, switching window size $L$, and tolerance $\epsilon_{\text{ROMSAD}}$.
\IF{$1 \leq l \leq L$ \AND $\|\brho^{(l,*)} - \brho^{(l-1)}\| \geq \epsilon_{\text{ROMSAD}}$}
    \STATE Use ROMSA.
\ELSE
    \STATE Use DSA.
\ENDIF
\end{algorithmic}
\end{algorithm}



\subsection{Limitation of ROMSAD preconditioner \label{sec:limitation}}
The right-hand side of the ideal correction equation \eqref{eq:ideal_density_equation}, $\BSigma_s(\brho^{(l,*)}-\brho^{(l-1)})$, depends on the iteration number $l$. Hence, it implicitly relies on the initial guess and the underlying preconditioner. However, this preconditioner dependence is neglected by ROMSAD when building ROMs. This neglected dependence further leads to a mismatch between the offline and online residual trajectories. As a result, ROMSAD may suffer an online efficiency reduction after the first iteration.

\begin{figure}[!ht]
  \begin{center} 
\includegraphics[width=0.5\textwidth]{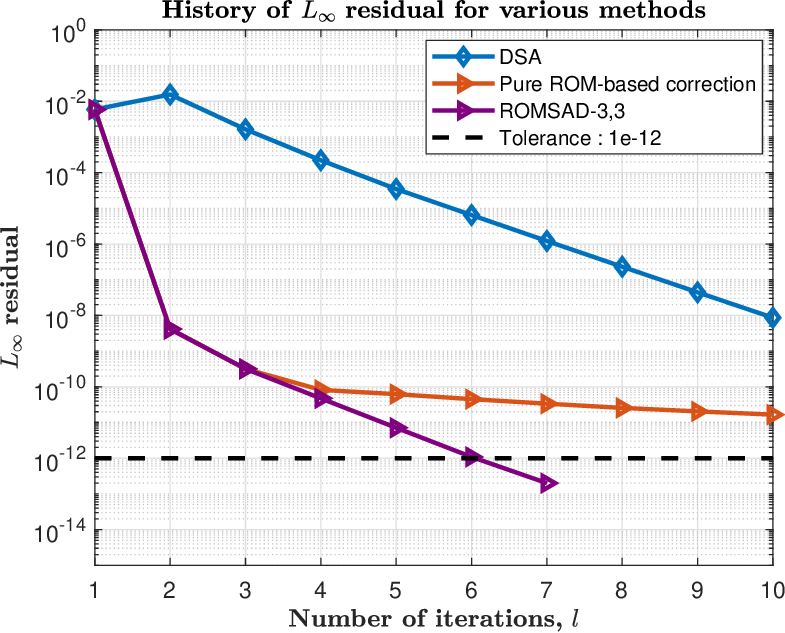}
  \caption{The residual history for the 1D setup in \eqref{eq:limitation_example} and $\mu = 39.0147$ with various correction strategies. Here, ROMSAD is implemented with SI. ROMSAD-$\mathfrak{w}, N_{w}$ denotes switching to DSA after the $N_{w}$-th iteration, with a snapshot window size of $\mathfrak{w}$. \label{fig:limitation}}
  \end{center}
\end{figure}

For illustration, we consider the following 1D parametric problem with zero inflow boundary conditions:
\begin{equation}
G(x)=\left\{
\begin{array}{ll}
1, & |x-1|\leq 0.5, \\
0, & \text{otherwise},
\end{array}
\right. \quad 
\sigma_{a}(x)=0, \quad 
\sigma_{s}(x)=0.1+\mu|x-1|, \quad 
x \in[0,2],
\label{eq:limitation_example}
\end{equation} 
where the parameter $\mu\in[10,100]$ controls how quickly the scattering grows from the domain center to the boundary. The ROM is built using $41$ uniformly sampled parameters from \ning{$[10,100]$}, with an SVD truncation threshold of $\epsilon_{\text{POD}}=10^{-9}$. As shown in Fig.~\ref{fig:limitation}, for a parameter outside the training set, using pure ROM-based corrections becomes progressively slower to converge. Leveraging the power of both ROM-based corrections and DSA, ROMSAD outperforms both of them. However, some efficiency reduction remains in ROMSAD after the first iteration. 

\begin{figure}[!ht]
    \centering
    \begin{subfigure}[b]{0.3\textwidth}
        \includegraphics[width=\textwidth]{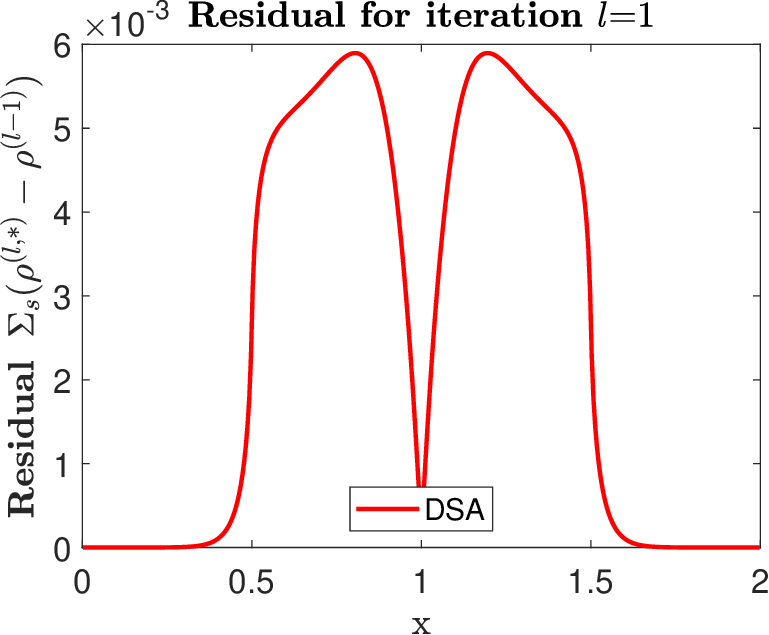} 
        \label{fig:dsa_k1}
    \end{subfigure}
    \hfill
    \begin{subfigure}[b]{0.3\textwidth}
        \includegraphics[width=\textwidth]{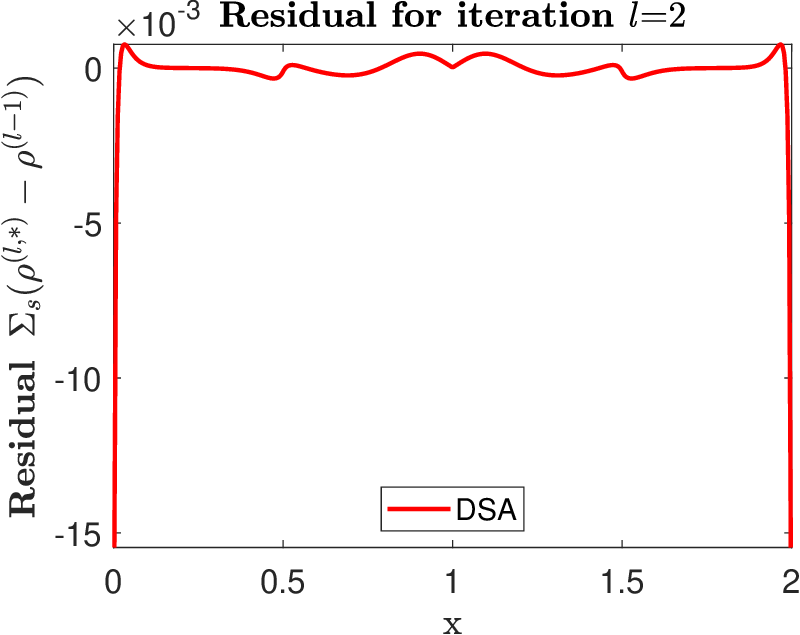} 
        \label{fig:dsa_k2}
    \end{subfigure}
    \hfill
    \begin{subfigure}[b]{0.3\textwidth}
        \includegraphics[width=\textwidth]{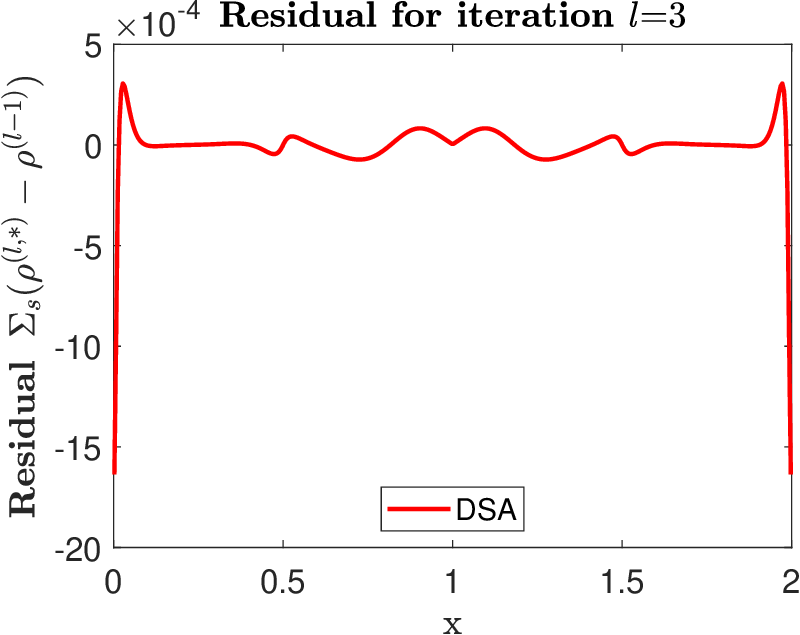} 
        \label{fig:dsa_k3}
    \end{subfigure}
    
    \begin{subfigure}[b]{0.3\textwidth}
        \includegraphics[width=\textwidth]{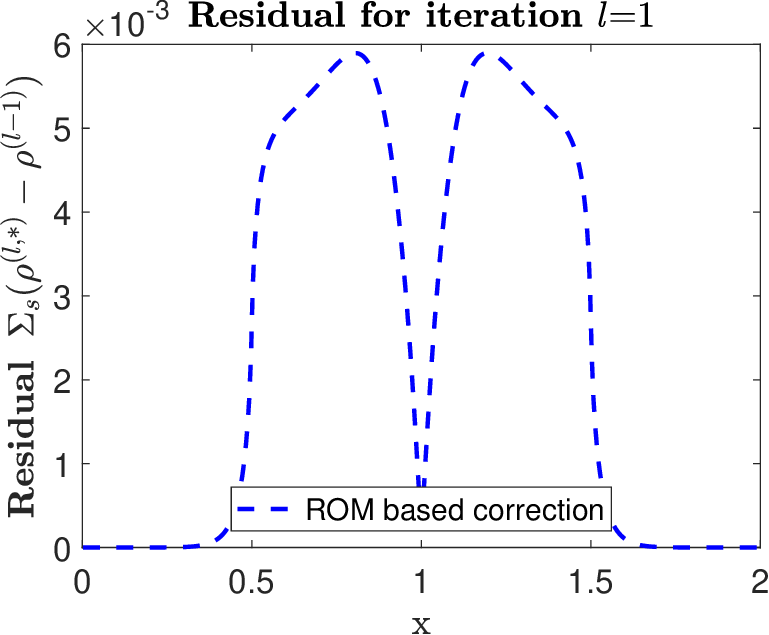} 
        \label{fig:rom_k1}
    \end{subfigure}
    \hfill
    \begin{subfigure}[b]{0.3\textwidth}
        \includegraphics[width=\textwidth]{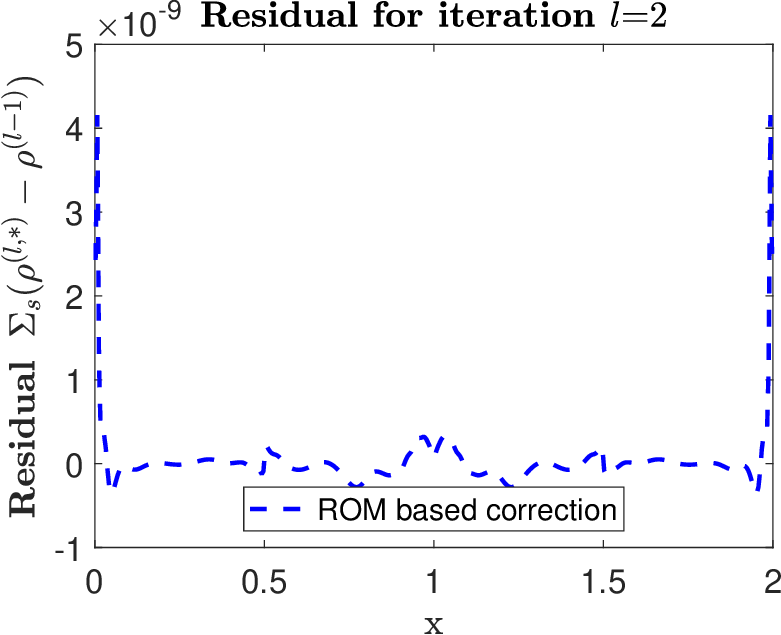}
        \label{fig:rom_k2}
    \end{subfigure}
    \hfill
    \begin{subfigure}[b]{0.3\textwidth}
        \includegraphics[width=\textwidth]{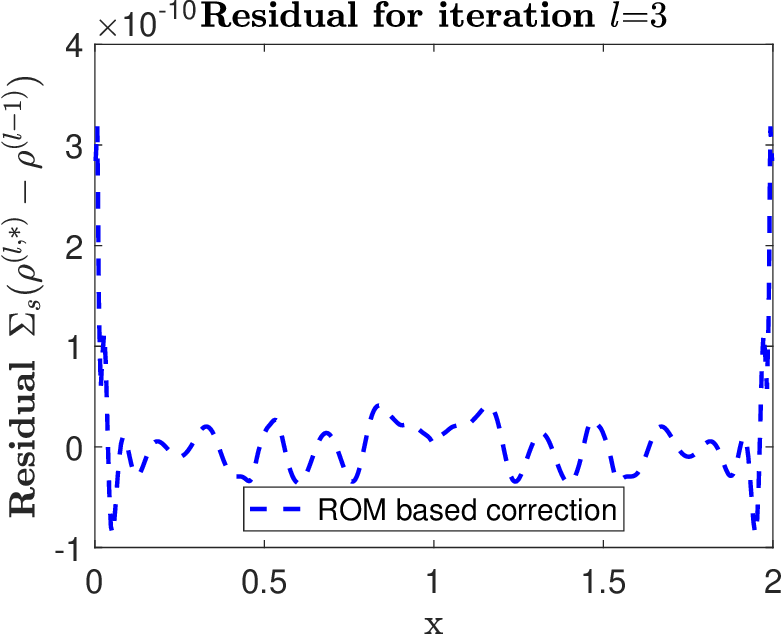} 
        \label{fig:rom_k3}
    \end{subfigure}
    
    \caption{The residual of \eqref{eq:source_iteration} in the SI step, $\BSigma_{s}(\brho^{(l,*)} - \brho^{(l-1)})$, for the 1D setup in \eqref{eq:limitation_example} and $\mu = 39.0147$ with DSA and preliminary ROM-based correction in \cite{peng2024romsad}. \label{fig:residual_difference}}
\end{figure}

This efficiency loss arises from the offline ROM construction. In particular, correction snapshots are generated using a fixed preconditioner (e.g., DSA) offline. As a result, the residual trajectory encoded in the ROM reflects the offline preconditioner, not the online residuals produced by the ROM-based preconditioner itself. In Fig.~\ref{fig:residual_difference}, this offline-online mismatch is visually demonstrated: while the residuals match in the first iteration, substantial deviations in both shape and magnitude exist in later iterations. This observation explains why ROMSAD performs well initially but suffers efficiency degradation in subsequent iterations.

\begin{rem}
Here, the term $\BSigma_{s}(\brho^{(l,*)} - \brho^{(l-1)})$ represents the residual for the SI equation in \eqref{eq:dg_matrix_vec} at the $\bff$ level, while $\brho^{(l,*)} - \brho^{(l-1)}$ is the residual at the $\brho$ level, as derived in \eqref{eq:rho_residual}.
\end{rem}

\section{Trajectory-Aware ROM framework for SI-SA \label{sec:trajectory_propose}}
According to Sec.~\ref{sec:limitation}, ROMSAD overlooks the preconditioner dependence of the residual trajectory when building ROMs for the ideal correction equation \eqref{eq:ideal_density_equation}, causing a significant mismatch between the online and offline residual trajectories after the first iteration. To mitigate this offline-online trajectory mismatch, we develop a trajectory-aware framework to iteratively construct ROMs \ning{by adapting the multispace reduced basis (MSRB) preconditioning \cite{santo2018multi} to parametric RTEs.}
In this section, we first present the SI-based implementation of the proposed framework \ning{(SI-TAR)}, then extend it to FGMRES \ning{(FGMRES-TAR)}, discuss practical implementation details, \ning{and clarify the connections and differences between our method and MSRB \cite{santo2018multi}.}

\subsection{Trajectory-Aware ROM for SI}
We first outline the trajectory-aware ROM framework and then present its main components for the SI solver (SI-TAR). The full offline and online algorithms are summarized in Alg.~\ref{alg:tar-romig-1} and Alg.~\ref{alg:tar-romig-2}, respectively.

\begin{figure}[!ht]
\centering

\begin{subfigure}[b]{\textwidth}
\centering
\begin{tikzpicture}[
    >=Stealth,
    node distance=0.8cm and 0.8cm,
    every node/.style={font=\footnotesize}
]
\node (p0) {$\brho^{(0)}$};
\node [right=of p0] (p1) {$\bff^{(1)},\; \brho^{(1,*)}$};
\node [right=of p1] (p2) {\fbox{$\delta\brho^{(1)}=\textcolor{red}{\BC^{-1}_{\mathrm{DSA}}}\BSigma_{s}(\brho^{(1,*)}-\brho^{(0)})$}};
\node [right=of p2] (p3) {$\brho^{(1)}_{\mathrm{DSA}}$};

\node [below=0.4cm of p3] (p4) {$\bff_{\mathrm{DSA}}^{(2)},\;\brho_{\mathrm{DSA}}^{(2,*)}$};
\node [left=of p4] (p5) {\fbox{$\delta\brho^{(2)}_{\mathrm{DSA}}=\textcolor{red}{\BC^{-1}_{\mathrm{DSA}}}\BSigma_{s}(\brho^{(2,*)}_{\mathrm{DSA}}-\brho^{(1)}_{\mathrm{DSA}})$}};
\node [left=0.5cm of p5] (p6) {$\cdots$};
\node [left=1.2cm of p6] (p7) {\textcolor{blue}{$\BC^{-1}_{\mathrm{ROM}}$}};

\draw[->] (p0) -- node[above,font=\scriptsize] {SI} (p1);
\draw[->] (p1) -- (p2);
\draw[->] (p2) -- (p3);

\draw[->] (p4) -- (p5);
\draw[->] (p5) -- (p6);
\draw[->] (p6) -- node[above,font=\scriptsize] {Output}(p7);

\draw[->] (p3) -- node[right,font=\scriptsize] {SI} (p4);

\draw[gray!50, dashed, thick] 
(current bounding box.south west) rectangle 
(current bounding box.north east);

\end{tikzpicture}
\caption{ROM-enhanced preconditioner hybridized with DSA (ROMSAD)}
\end{subfigure}
\hfill

\begin{subfigure}[b]{\textwidth}
\centering
\begin{tikzpicture}[
    >=Stealth,
    node distance=0.8cm and 0.8cm,
    every node/.style={font=\footnotesize}
]

\node (p0) {$\brho^{(0)}$};
\node [right=of p0] (p1) {$\bff^{(1)},\; \brho^{(1,*)}$};
\node [right=2cm of p1] (p2) {\fbox{$\delta\brho^{(1)}_{\mathrm{ROM}}=\textcolor{blue}{\BC^{-1}_{\mathrm{ROM,1}}}\BSigma_{s}(\brho^{(1,*)}-\brho^{(0)})$}};
\node [right=of p2] (p3) {$\brho_{\mathrm{ROM}}^{(1)}$};

\node [below=0.4cm of p3] (p4) {$\bff^{(2)}_{\mathrm{ROM}},\;\brho_{\mathrm{ROM}}^{(2,*)}$};
\node [left=1.8cm of p4] (p5) {\fbox{$\delta\brho^{(2)}_{\mathrm{ROM}}=\textcolor{blue}{\BC^{-1}_{\mathrm{ROM,2}}}\BSigma_{s}(\brho^{(2,*)}_{\mathrm{ROM}}-\brho^{(1)}_{\mathrm{ROM}})$}};
\node [left=0.5cm of p5] (p6) {$\cdots$};
\node [left=1.2cm of p6] (p7) {$\textcolor{blue}{\{\BC_{\mathrm{ROM},l}^{-1}\}_{l=1}^{N_{w}}}$};

\draw[->] (p0) -- node[above,font=\scriptsize] {SI} (p1);
\draw[->] (p1) --  node[above,font=\scriptsize] {Construct} node[below,font=\scriptsize] {$\BC^{-1}_{\mathrm{ROM,1}}$} (p2);
\draw[->] (p2) --  (p3);

\draw[->] (p4) -- node[above,font=\scriptsize] {\;\;Construct} node[below,font=\scriptsize] {\;\;$\BC^{-1}_{\mathrm{ROM,2}}$}(p5);
\draw[->] (p5) -- (p6);
\draw[->] (p6) -- node[above,font=\scriptsize] {Output} (p7);

\draw[->] (p3) --
node[right,font=\scriptsize] {SI}  (p4);

\draw[gray!50, dashed, thick] 
(current bounding box.south west) rectangle 
(current bounding box.north east);

\end{tikzpicture}
\caption{Trajectory-aware ROM (TAR)}
\end{subfigure}
\caption{The flowcharts of the offline stage for the ROMSAD method and the trajectory-aware method. Here, $\BC_{\textrm{ROM},l}^{-1}$ and $\BC_{\textrm{DSA}}^{-1}$ are operators for the ROM-based and DSA-based density correction.\label{fig:flow}}
\end{figure}
\textbf{Key ideas and motivation:}
In Sec.~\ref{sec:limitation}, we demonstrate that the reduced online efficiency of ROMSAD arises from a mismatch between the offline residual trajectory generated by DSA and the online residual trajectory determined by the ROM itself. 

\ning{To mitigate this mismatch, we adapt MSRB preconditioning to parametric RTEs by sequentially constructing iteration-dependent ROMs that account for the influence of preceding ROM corrections on subsequent residual trajectories.} This construction aims to better align the offline and online residual trajectories. Fig.~\ref{fig:flow} illustrates the proposed offline strategy in comparison with ROMSAD. We emphasize that the essential advance is not simply using different ROMs across iterations, but their sequential and trajectory-aware construction to ensure the consistency between the offline and online residual trajectories.

\textbf{Trajectory-aware offline stage:}
For the SI solver, we sequentially construct iteration-dependent ROMs within the total aware level $N_w$, rather than using a fixed SI-DSA trajectory to generate correction snapshots. Each ROM is constructed along the residual trajectory influenced by the ROM-based corrections applied in preceding iterations. 

Specifically, we first save converged solutions for a training set, namely $\bff_{\bmu_1},\dots,\bff_{\bmu_{N_{\textrm{train}}}}$. Then, given an initial guess $\brho_{\bmu}^{(0)}$, we perform one SI step solving
\begin{equation}
(\BD_j+\BSigma_t)\bff_{\bmu,j}^{(1)}=\BSigma_s\brho_{\bmu}^{(0)}+\widetilde{\BQ}_{\bmu,j},\;j=1,\dots,N_{\Bupsilon},
\end{equation}
and  construct a correction snapshot matrix $\BF_{c}^{(1)}\in \mathbb{R}^{N_{h}\times N_{\textrm{train}}}$ defined as:
\begin{equation}
    \BF_{c}^{(1)} = \left[\delta\bff^{(1)}_{\bmu_1}, \delta\bff^{(1)}_{\bmu_2},\cdots,\delta\bff^{(1)}_{\bmu_{N_{\mathrm{train}}}}\right],\;\;\text{where}\;\;\delta\bff^{(1)}_{\bmu_i} = \bff_{\bmu_{i}} -  \bff_{\bmu_{i}}^{(1)}.    \label{eq:correction_snapshot}
\end{equation}
Next, we apply SVD to the snapshot matrix $\BF_{c}^{(1)}$ to construct the reduced basis for the first iteration $\BU_{r_{1}}^{(1)}\in\mathbb{R}^{\ning{N_h\times r_1}}$. Before the next iteration, we introduce a ROM-based  correction: $\brho^{(1)}_{\bmu}=\brho^{(1,*)}_{\bmu}+\BU_{r_{1}}^{\rho,(1)}\delta\bc_{\bmu}^{(1)}$, where 
\begin{equation}
(\BU_{r_{1}}^{(1)})^T\BA_{\bmu}\BU_{r_{1}}^{(1)}\delta\bc_{\bmu}^{(1)}=(\BU_{r_{1}}^{(1)})^T\delta\bb_{\bmu}^{(1)},
\label{eq:tar-rom-correction}
\end{equation}
and $\BU_{r_{1}}^{\rho,(1)}=\sum_{j=1}^{N_{\Bupsilon}}\omega_j\BU_{r_1,j}^{(1)}$ with $\BU_{r_{1},j}^{(1)}$ defined similarly to Sec.~\ref{sec:rom}. 
Following the same methodology, this procedure can be iteratively applied for subsequent iterations $(l\geq2)$ until reaching the maximum aware level $N_w$. This iterative offline process provides a sequence of reduced bases $\{\BU_{r_{l}}^{(l)}\}_{l=1}^{N_w}$ with $\BU_{r_{l}}^{(l)}\in\mathbb{R}^{\ning{N_h\times r_l }}$ and corresponding reduced-order operators. 

The key distinction from ROMSAD lies in the sequential construction of ROMs using preceding ROM-based corrections \eqref{eq:tar-rom-correction}. This strategy mitigates the offline–online trajectory mismatch.

\textbf{Trajectory-aware online stage:} For a test parameter $\bmu\in\mathscr{P}_{\mathrm{test}}$, we first perform a transport sweep to obtain $\brho_{\bmu}^{(l,*)}$ from $\brho_{\bmu}^{(l-1)}$. During the first $N_w$ iterations, the density is then updated using the corresponding ROM-based correction:
$\brho_{\bmu}^{(l)} = \brho_{\bmu}^{(l,*)} +\BU_{r_l}^{\rho,(l)}\delta\bc_{\bmu}^{(l)}$, where
\begin{equation}
(\BU_{r_l}^{(l)})^T
\BA_{\bmu}
\BU_{r_l}^{(l)}
\delta\bc_{\bmu}^{(l)}
=
(\BU_{r_l}^{(l)})^T
\delta\bb_{\bmu}^{(l)}.
\end{equation}
Starting from the $(N_w+1)$-th iteration, the ROM-based correction is replaced by DSA.

\textbf{Initial guess:} Our trajectory-aware framework allows for different initial guesses. In addition to a zero initial
guess, we consider the scalar-flux approximation obtained from the reduced
system \eqref{eq:reduced_equation}. This ROM-based initial guess is inexpensive
to compute and can reduce the number of online iterations. 

\begin{algorithm}[ht]
\caption{Offline stage of trajectory-aware ROMs for SI (\ning{SI-TAR}). \label{alg:tar-romig-1} }
\begin{algorithmic}[1]
\STATE{\underline{\textit{\textbf{Offline stage:}}}}
\STATE{Given an initial guess $\brho^{(0)}_{\mathrm{train}}$, converged solutions namely $\ning{\bff}_{\bmu}$, and the total aware level $N_w$.}
\FOR{ $l = 1:N_{w}$}
    \STATE{\textbf{Perform source iteration:} $\forall\bmu\in\mathscr{P}_{\mathrm{train}}$,  get $\bff^{(l)}_{\bmu}$ by solving  
    \begin{equation}
       (\BD_j+\BSigma_{t,\bmu}) \bff_{\bmu,j}^{(l)} =\BSigma_{s,\bmu}\brho^{(l-1)}_{\bmu}+\widetilde{\bQ}_{\bmu,j},\;\; j=1,\dots,N_{\Bupsilon}, \;\;\forall\;l\geq 1, 
    \end{equation}
    through transport sweeps with given $\brho_{\bmu}^{(l-1)}$, and calculate $\brho^{(l,*)}_{\bmu}$. }
    \STATE{\textbf{Construct \textit{trajectory-aware} correction snapshots:}  $\BF_{c}^{(l)}=[\delta\bff_{\bmu}^{(l)}]$ as in equation \eqref{eq:correction_snapshot}.}
    \STATE{\textbf{Construct the reduced basis:} apply SVD to the snapshot $\BF_{c}^{(l)}$ to get the reduced basis $\BU_{r_{l}}^{(l)}$.}
    \STATE{\textbf{Attain \textit{trajectory-aware} density correction:}  $\forall\bmu\in\mathscr{P}_{\mathrm{train}}$, solve the reduced-order system
    \begin{equation}(\BU_{r_{l}}^{(l)})^{T}\BA_{\bmu}\BU^{(l)}_{r_{l}}\delta\mathbf{c}_{\bmu,r_{l}}^{(l)} = (\BU^{(l)}_{r_{l}})^{T}\delta\mathbf{b}_{\bmu}^{(l)}.
     \end{equation}
    }
    \STATE{Update density as $\brho^{(l)}_{\bmu} = \brho^{(l,*)}_{\bmu} + \delta\brho^{(l)}_{\bmu}$, where
    \begin{equation}
    \delta\brho^{(l)}_{\bmu}=\sum_{j=1}^{N_{\Bupsilon}}\omega_j\delta \bff_{\ning{\bmu,j}}^{(l)}=  \sum_{j=1}^{N_{\Bupsilon}}
\omega_j\BU^{(l)}_{r_l,j}\delta\bc^{(l)}_{\bmu,r_l} = \mathbf{U}_{r_{l}}^{\rho,(l)}\delta\bc^{(l)}_{\bmu,r_l}.
\end{equation}}
\ENDFOR
\STATE{Return the reduced bases $\{\BU_{r_{l}}^{(l)}\}_{l=1}^{N_w}$ and $\{\BU_{r_{l}}^{\rho,(l)}\}_{l=1}^{N_w}$ for related reduced-order operators.}
\end{algorithmic}
\end{algorithm}

\begin{algorithm}[ht]
\caption{Online stage of trajectory-aware ROMs for SI (\ning{SI-TAR}). \label{alg:tar-romig-2} }
\begin{algorithmic}[1]
\STATE{\underline{\textit{\textbf{Online stage:}}}}
\STATE{Given an initial guess $\brho^{(0)}_{\mathrm{test}}$, total aware level $N_w$, tolerance $\epsilon_{\mathrm{tar}}$ and maximum iteration number $N_k$.}
\FOR{$ k = 1: N_k$}
    \STATE{\textbf{Perform source iteration:} $\forall\hat{\bmu}\in\mathscr{P}_{\mathrm{test}}$, get $\brho^{(k,*)}_{\hat{\bmu}}$ and $\bff^{(k)}_{\hat{\bmu}}$ with input $\brho^{(k-1)}_{\hat{\bmu}}$.}
    \IF{ $\|\brho^{(k,*)}_{\hat{\bmu}} - \brho^{(k-1)}_{\hat{\bmu}}\|_{2}/\|\tilde{\bb}_{\hat{\bmu}}\|_2 < \epsilon_{\mathrm{tar}}$}
        \STATE{Return solutions: $\brho^{(k,*)}_{\hat{\bmu}}$ and $\bff^{(k)}_{\hat{\bmu}}$}.
    \ENDIF
    \STATE{\textbf{Select and solve a surrogate model:}}
    \IF{$k\leq N_w$}
        \STATE{Apply \textit{trajectory-aware} ROMs: \begin{equation}(\BU_{r_{k}}^{(k)})^{T}\BA_{\hat{\bmu}}\BU^{(k)}_{r_{k}}\delta\mathbf{c}_{\hat{\bmu},r_{k}}^{(k)} = (\BU^{(k)}_{r_{k}})^{T}\delta\mathbf{b}_{\hat{\bmu}}^{(k)},
        \label{eq:tar_ig_rom}\end{equation} 
        and attain the correction $\delta\rho_{\hat{\bmu}}^{(k)}=\BU_{r_k}^{\brho,(k)}\delta\bc_{\hat{\bmu},r_k}^{(k)}$.}
    \ELSE
        \STATE{Apply DSA to get the correction $\delta\brho_{\hat{\bmu}}^{(k)}$.}
    \ENDIF
    \STATE{Update density as $\brho^{(k)}_{\hat{\bmu}} = \brho^{(k,*)}_{\hat{\bmu}} + \delta\brho^{(k)}_{\hat{\bmu}}$.}
\ENDFOR
\end{algorithmic}
\end{algorithm}

\begin{rem}\label{rem:offline-cost}
Compared to the ROMSAD method in \cite{peng2024romsad}, the offline stage of the trajectory-aware approach requires additional $N_w$ transport sweeps for each training parameter. Fortunately, we find that using $N_w=1,2$ already leads to significant online acceleration over DSA and ROMSAD in practice.  
\end{rem}
\vspace{-15pt}
\ning{\begin{rem}[Conditioning of the Galerkin reduced systems] Since $A_\mu$ is non-symmetric, the Galerkin reduced matrices $(U_r^{(l)})^T A_\mu U_r^{(l)}$ are not guaranteed to be well conditioned. In \ref{sec:condition_number}, we report the condition numbers of the Galerkin reduced systems for the variable-scattering test (Sec.~\ref{sec:scattering}), which remain $\mathcal{O}(10^3)$ at most across all test parameters, indicating well-behaved conditioning. To improve robustness, least-squares Petrov-Galerkin projection \cite{behne2021model,matsuda2025reduced} can be considered as an alternative.
\end{rem}}

\subsection{Trajectory-Aware ROM for flexible GMRES}
As discussed in \cite{warsa2004krylov,peng2025flexible}, Krylov methods generally exhibit more robust convergence than SI in transport problems. In what follows, we extend the trajectory-aware framework to the flexible generalized minimal residual (FGMRES) method, yielding FGMRES-TAR.

\subsubsection{Ideal correction within Krylov framework\label{sec:ideal_correction_krylov}}
FGMRES \cite{saad1993flexible} allows iteration-dependent preconditioners, as summarized in Alg.~\ref{alg:fgmres}.  To extend our trajectory-aware ROMs to FGMRES,  we first formulate the  ideal correction equation for FGMRES and then show how to obtain its solution without directly solving the global correction system.

\textbf{Ideal correction equation:} As shown in \eqref{eq:ideal_density_equation} and \eqref{eq:rho_residual}, the ideal correction equation for SI is 
\begin{equation}
    (\BD_j+\BSigma_t)\delta\bff_j^{(l)} = \BSigma_{s}\delta\brho^{(l)}+\BSigma_s(\brho^{(l,*)}-\brho^{(l-1)})= \BSigma_{s}\delta\brho^{(l)}+\BSigma_s\br^{(l-1)}.
\end{equation}
The SA preconditioner acts on the residual $\br^{(l-1)}$ in SI, while it acts on the Krylov vector $\bq^{(l)}$ in FGMRES. Accordingly, the ideal correction for FGMRES becomes 
\begin{equation}
    (\BD_j+\BSigma_t)\delta\bff_j^{(l)} = \BSigma_{s}\delta\brho^{(l)} +  \BSigma_s\bq^{(l)}.
    \label{eq:ideal_density_equation_fgmres}
\end{equation}
Following \cite{peng2025flexible} (see \ref{sec:equation_equivalence} for details), the ideal correction equation \eqref{eq:ideal_density_equation_fgmres} is equivalent to
\begin{equation}
    (\BD_j+\BSigma_t) \delta\bff_j^{(l)} =\BSigma_s\beeta^{(l)},\;\; j=1,\dots,N_{\Bupsilon}, 
    \label{eq:equivalence_1}
\end{equation}
where $\beeta^{(l)}$ solves the following  equation:
\begin{equation}
    \WBA\beeta^{(l)} = (\BI-\mathbf{K}\BSigma_{s})\beeta^{(l)} = \bq^{(l)}.
    \label{eq:equivalence_2}
\end{equation}
 
\textbf{Obtain the ideal correction without solving the correction equation:} Similar to \cite{santo2018multi,peng2025flexible}, the auxiliary vector $\beeta^{(l)}$ in \eqref{eq:equivalence_1} can be computed without solving \eqref{eq:equivalence_2}. For the first iteration, we have 
\begin{subequations}
    \label{eq:eta1}
    \begin{align}
    \WBA\beeta^{(1)} &= (\BI-\mathbf{K}\BSigma_{s})\beeta^{(1)} = \bq^{(1)} = \br^{(0)}/ \|\br^{(0)}\|_{2} = (\wbb-\WBA\brho^{(0)})/\|\br^{(0)}\|_{2},\\
    \beeta^{(1)} &= \WBA^{-1}(\wbb-\WBA\brho^{(0)})/\|\br^{(0)}\|_{2} = \WBA^{-1}(\WBA\brho-\WBA\brho^{(0)})/\|\br^{(0)}\|_{2} = (\brho-\brho^{(0)})/\|\br^{(0)}\|_{2},
    \end{align}
\end{subequations}
where $\brho$ is the converged solution to $\WBA\brho=\wbb$ given by FGMRES. In later iterations ($l\geq2$), the Arnoldi relation yields the following recurrence for $\beeta^{(l)}$: 
\begin{subequations}
    \label{eq:eta2}
    \begin{align}
    &\WBA\beeta^{(l)} = \bq^{(l)} = \bw/ \mathbf{H}_{l,l-1} = (\WBA \bz^{(l-1)} - \sum_{k=1}^{l-1}\mathbf{H}_{k,l-1}\bq^{(k)})/ \mathbf{H}_{l,l-1},\quad l\geq 2,\\
    & \beeta^{(l)} = (\bz^{(l-1)} - \sum_{k=1}^{l-1}\mathbf{H}_{k,l-1}\WBA^{-1}\bq^{(k)})/ \mathbf{H}_{l,l-1} = (\bz^{(l-1)} - \sum_{k=1}^{l-1}\mathbf{H}_{k,l-1}\beeta^{(k)})/ \mathbf{H}_{l,l-1},\quad l\geq 2.
    \end{align}
\end{subequations} 
Following the above derivation, the ideal corrections $\delta\bff^{(l)}$ can be computed as follows:
\begin{enumerate}
    \item Store the converged solution $\brho$, the vectors $\bz^{(l)}$, and the \ning{Hessenberg} matrix entries for $1\leq l\leq N_w$. Then, recursively compute $\beeta^{(l)}$ using \eqref{eq:eta1} and \eqref{eq:eta2}.
    \item Given $\beeta^{(l)}$, compute $\delta\bff^{(l)}$ by solving \eqref{eq:equivalence_1} via transport sweeps.
\end{enumerate}

\begin{algorithm}[!ht]
\caption{Flexible GMRES method to solve $\WBA\brho=\wbb$ \eqref{eq:matrix_vec_rewrite_rho} with a preconditioner $\BM_{j}^{-1}$ \cite{saad1993flexible}. \label{alg:fgmres} }
\begin{algorithmic}[1]
\STATE{Given an initial guess $\brho^{(0)}$, a preconditioner $\BM^{-1}_{j}$, and the maximum number of iterations $n$.}
\STATE{\textbf{Initialization:} $\br^{(0)}=\wbb-\WBA\brho^{(0)}$.} 

\STATE{Define $\beta=||\br^{(0)}||_2$ and $\boldsymbol{q}^{(1)}=\br^{(0)}/\beta$. Allocate the memory for the Hessenberg matrix $\mathbf{H}\in\mathbb{R}^{(n+1)\times n}$.}
\STATE{\textbf{\ning{Arnoldi} process:}}
\FOR{$l=1:n$}
\STATE{\textbf{Apply the \ning{iteration-dependent} preconditioner:}  $\bz^{(l)}:=\BM_{l}^{-1}\mathbf{q}^{(l)}$.}
\STATE{Compute $\bw:=\WBA\bz^{(l)}$.}
    \FOR{$i=1:l$}
        \STATE{\;\;$\mathbf{H}_{il}=\bw^T \mathbf{q}^{(i)}$, and update 
        $\bw=\bw-\mathbf{H}_{il}\mathbf{q}^{(i)}$.}
    \ENDFOR
\STATE{Compute $\mathbf{H}_{l+1,l}=||\bw||_2$.}
\STATE{Update $\mathbf{q}^{(\ning{l}+1)}=\bw/\mathbf{H}_{l+1,l}$.}
\IF{the stopping criterion is satisfied}
    \STATE{Break}
\ENDIF
\ENDFOR
\STATE{\textbf{Obtain the solution:} Define $\widetilde{\BZ}=\left(\bz^{(1)},\bz^{(2)},\cdots,\bz^{(l)}\right)$, $\widetilde{\mathbf{H}} = \left(\mathbf{H}^{(1)},\mathbf{H}^{(2)},\cdots,\mathbf{H}^{(l)}\right)$ and solve the minimization problem:
\begin{equation}
    \by^{*} = \arg\min_{\by\in \mathbb{R}^{l}}\|\beta\be_{1}-\widetilde{\mathbf{H}}\by\|_2, \quad \text{where}\;\; \be_{1} = (1,0,\cdots,0),
\end{equation}
then return the solution $\brho:=\brho^{(0)}+\widetilde{\BZ}\by^*$} 
\end{algorithmic}
\end{algorithm}
\subsubsection{ROM-enhanced SA preconditioner}
Following Sec.~\ref{sec:sa-krylov}, we derive the SA preconditioner for the Krylov method based on the correction strategy of SI. 
Assume the reduced basis for the $l$-th iteration is $\BU_{r_{l}}^{(l)}\in\mathbb{R}^{N_{\Bupsilon}N_{\mathrm{DOF}}\times \ning{r_{l}}}$, where $r_{l}$ denotes the dimension of the reduced basis at the $l$-th iteration.  For simplicity, we also define $\BU_{r_{l},j}^{(l)}\in \mathbb{R}^{N_{\mathrm{DOF}}\times r_{l}}$ as the row submatrix of $\BU_{r_{l}}^{(l)}$.
Then, the ROM-based correction for SI is as follows. 
\begin{subequations}
\begin{align}
&(\BU_{r_{l}}^{(l)})^{T}\BA_{\bmu}\BU^{(l)}_{r_{l}}\delta\mathbf{c}_{\bmu,r_{l}}^{(l)} = (\BU^{(l)}_{r_{l}})^{T}\delta\mathbf{b}_{\bmu}^{(l)} = \sum_{j=1}^{N_{\Bupsilon}} (\BU^{(l)}_{r_{l},j})^{T}\BSigma_{s}(\brho^{(l,*)} - \brho^{(l-1)}) = (\sum_{j=1}^{N_{\Bupsilon}} \BU^{(l)}_{r_{l},j})^{T}\BSigma_{s}\br^{(l-1)},\\
& \delta\bff^{(l)} \approx \BU^{(l)}_{r_{l}}\delta\bc_{\bmu,r_{l}}^{(l)},\quad \delta\brho^{(l)} \approx \sum_{j=1}^{N_{\Bupsilon}} \omega_{j}\BU^{(l)}_{r_{l},j}\delta\bc_{\bmu,r_{l}}^{(l)} = (\sum_{j=1}^{N_{\Bupsilon}} \omega_{j}\BU^{(l)}_{r_{l},j}) \delta\bc_{\bmu,r_{l}}^{(l)},
\end{align}
\label{eq:tar_cor_rho}
\end{subequations}
Following \eqref{eq:tar_cor_rho}, we have
\begin{equation}
    \delta\brho^{(l)} \approx (\sum_{j=1}^{N_{\Bupsilon}} \omega_{j}\BU^{(l)}_{r_{l},j}) \delta\bc_{\bmu,r_{l}}^{(l)} =  (  \sum_{j=1}^{N_{\Bupsilon}} \omega_{j}\BU^{(l)}_{r_{l},j})\left((\BU_{r_{l}}^{(l)})^{T}\BA_{\bmu}\BU^{(l)}_{r_{l}}\right)^{-1}(\sum_{j=1}^{N_{\Bupsilon}} \BU^{(l)}_{r_{l},j})^{T}\BSigma_{s}\br^{(l-1)}.
    \label{eq: delta_rho_tar}
\end{equation}
By comparing \eqref{eq: delta_rho_tar} and \eqref{eq:delta_rho_c}, the reduced-order operator for the $l$-th iteration is defined as:
\begin{equation}
    \BC_{\mathrm{ROM}, l}^{-1} = (  \sum_{j=1}^{N_{\Bupsilon}} \omega_{j}\BU^{(l)}_{r_{l},j})\left((\BU_{r_{l}}^{(l)})^{T}\BA_{\bmu}\BU^{(l)}_{r_{l}}\right)^{-1}(\sum_{j=1}^{N_{\Bupsilon}} \BU^{(l)}_{r_{l},j})^{T},
\end{equation}
Similar to \eqref{eq:left_preconditioner}, the corresponding ROM-enhanced SA preconditioner is then given by:
\begin{equation}
    \BM_{\mathrm{ROM},l}^{-1} = \BI + \BC_{\mathrm{ROM}, l}^{-1}\BSigma_{s}  = \BI + (  \sum_{j=1}^{N_{\Bupsilon}} \omega_{j}\BU^{(l)}_{r_{l},j})\left((\BU_{r_{l}}^{(l)})^{T}\BA_{\bmu}\BU^{(l)}_{r_{l}}\right)^{-1}(\sum_{j=1}^{N_{\Bupsilon}} \BU^{(l)}_{r_{l},j})^{T}\BSigma_{s}.
    \label{eq:preconditioner_M}
\end{equation}

Using the same switching strategy as the SI case, we utilize ROM-based preconditioners in the first few iterations, followed by a  transition to DSA to ensure robustness: 
\begin{equation}
    \BM_{l}^{-1} = 
\begin{cases} 
\BM_{\mathrm{ROM},l}^{-1}, & 1 \leq l \leq N_{w}, \\
\mathbf{M}_{\mathrm{DSA}}^{-1}, & \text{otherwise.}
\label{eq:tar_preconditioner}
\end{cases}
\end{equation}
Unlike in our previous work \cite{peng2025flexible}, $\BM_{\mathrm{ROM},l}^{-1}$ varies in every iteration with $1\leq l\leq N_w$.

\subsubsection{Offline and online stages of FGMRES-TAR \label{sec:fgmres-offline}}
\textbf{Offline stage:} When extending the trajectory-aware approach to FGMRES, we follow the same principle as in SI-TAR. The key is not merely to use different ROMs across iterations, but to construct them sequentially based on preceding ROM-based corrections. The main difference from SI-TAR lies in the construction of correction snapshots, as described in Sec.~\ref{sec:ideal_correction_krylov}. The complete offline algorithm is summarized in Alg.~\ref{alg:tar-gmres}.

\begin{algorithm}[!ht]
\caption{Offline stage of trajectory-aware ROMs for FGMRES solver (\ning{FGMRES-TAR}).\label{alg:tar-gmres}}
\begin{algorithmic}[1]
\STATE{\underline{\textit{\textbf{Offline stage:}}}}
\STATE Given an initial guess $\brho^{(0)}_{\mathrm{train}}$, the converged solutions $\brho_{\mathrm{train}}$, and total aware level $N_{w}$.
\STATE{\textbf{Initialization:} $\br^{(0)}=\wbb-\WBA\brho^{(0)}$.} 
\STATE Define $\beta=\left\|\mathbf{r}^{(0)}\right\|_2$ and $\mathbf{q}^{(1)}=\mathbf{r}^{(0)}/ \beta$. Allocate the memory for the Hessenberg matrix $\mathbf{H} \in \mathbb{R}^{(N_{w}+1) \times N_{w}}$.
\STATE \textbf{Preconditioner Generation:}
\FOR {$l = 1 : N_{w}$}
\STATE{\textbf{Construct \textit{trajectory-aware} correction snapshot:}$\;\forall \bmu\in\mathscr{P}_{\textrm{train}}$,}
\IF {$l = 1$}
\STATE $\beeta^{(l)}_{\bmu} = (\brho_{\bmu}-\brho^{(0)}_{\bmu})/\left\|\mathbf{r}^{(0)}\right\|_2.$
\ELSE
\STATE $\beeta^{(l)}_{\bmu} = \left(\mathbf{z}_{\bmu}^{(l-1)}-\sum_{k=1}^{l-1} \mathbf{H}_{k, l-1,\bmu} \beeta^{(k)}_{\bmu}\right)/\mathbf{H}_{l, l-1,\bmu}$.
\ENDIF
\STATE{ Solve $(\BD_{j}+\BSigma_{t,\bmu})\delta \bff_{j,\bmu}^{(l)}=\BSigma_{s,\bmu}\beeta_{\bmu}^{(l)}\;\;(j= 1,2,...,N_{v})$ through transport sweeps, and then assemble the snapshot matrix $\BF_{c}^{(l)}=[\delta\bff_{\bmu}^{(l)}]$ as in \eqref{eq:correction_snapshot}.}
    \STATE{\textbf{Construct the reduced basis:} apply SVD to the snapshot $\BF_{c}^{(l)}$ to get the reduced basis $\BU_{r_{l}}^{(l)}$, and update the preconditioner for the current iteration $\mathbf{M}^{-1}_{\mathrm{ROM},l}$ as in \eqref{eq:preconditioner_M}.}
\STATE \textbf{Apply the \ning{iteration-dependent} preconditioner:}  $\mathbf{z}_{\bmu}^{(l)}: = \mathbf{M}^{-1}_{\mathrm{ROM},l}\mathbf{q}_{\bmu}^{(l)}.$
\STATE \textbf{Continue the \ning{Arnoldi} process:}
\STATE Compute $\mathbf{w}_{\bmu}:=\widetilde{\mathbf{A}} \mathbf{z}_{\bmu}^{(l)}$.
\FOR {$k = 1 : l$}
\STATE $\mathbf{H}_{k,l,\bmu}=\mathbf{w}_{\bmu}^T \mathbf{q}_{\bmu}^{(\ning{k})}$, and update $\mathbf{w}_{\bmu}=\mathbf{w}_{\bmu}-\mathbf{H}_{k,l,\bmu} \mathbf{q}_{\bmu}^{(k)}$.
\ENDFOR
\STATE Compute $\mathbf{H}_{l+1, l,\bmu}=\|\mathbf{w}_{\bmu}\|_2$.
\STATE Update $\mathbf{q}_{\bmu}^{(l+1)}=\mathbf{w}_{\bmu} / \mathbf{H}_{l+1, l, \bmu}$. 
\IF {the stopping criterion is satisfied}
    \STATE{Break.}
\ENDIF
\ENDFOR
\STATE Return \ning{$\{\BU_{r_l}^{(l)}\}_{l=1}^{N_{w}}$} and store related quantities.
\end{algorithmic}
\end{algorithm}

\textbf{Online stage:} \ning{For each online parameter, we construct  $\BM_{\mathrm{ROM},l}^{-1}$ in \eqref{eq:preconditioner_M} using the stored reduced basis $\BU_{r_l}^{(l)}$. As outlined in Alg.~\ref{alg:fgmres}, we then apply $\BM_{\mathrm{ROM},l}^{-1}$ as the right preconditioner at iteration $l$, for $1\leq l\leq N_w$;} while for all subsequent iterations ($l>N_w$), we switch to the DSA preconditioner $\BM_{\mathrm{DSA}}^{-1}$. This hybrid approach combines the acceleration provided by ROM-based preconditioners with the robustness of DSA.

\ning{\subsection{Connections to and differences from MSRB}\label{sec:preconditioner_difference}}
 Our method \ning{is similar in spirit to} the \ning{multispace reduced basis (MSRB)} preconditioner for parametric elliptic and advection-diffusion equations in \cite{santo2018multi}. This preconditioner is also trajectory-aware. However, the two approaches differ in two important aspects.

First, our approach is derived from the SA formulation of the kinetic correction equation for RTEs, whereas MSRB \cite{santo2018multi} employs ROMs as coarse-level corrections within a general two-level preconditioning framework.
    
Second, compared to direct extensions of the multispace method in \cite{santo2018multi} to the \pzc{discrete integral scalar-flux formulation of} RTE \pzc{in equation \eqref{eq:rewrite_rho}}, our formulation is better suited to transport-sweep-based solvers for RTEs.
Within our trajectory-aware framework, the preconditioner for the $l$-th iteration \pzc{with $1\leq l\leq N_{w}$} takes the form:
    \begin{equation}
    \BM_{\mathrm{ROM},l}^{-1} = \BI + \BC_{\mathrm{ROM}, l}^{-1}\BSigma_{s}  = \BI + (  \sum_{j=1}^{N_{\Bupsilon}} \omega_{j}\BU^{(l)}_{r_{l},j})\left((\BU_{r_{l}}^{(l)})^{T}\BA_{\bmu}\BU^{(l)}_{r_{l}}\right)^{-1}(\sum_{j=1}^{N_{\Bupsilon}} \BU^{(l)}_{r_{l},j})^{T}\BSigma_{s}.
    \end{equation}
\ning{For comparison, we consider two direct extensions of \cite{santo2018multi} to the RTE with different fine-level preconditioners $\mathbf{P}_{\textrm{fine}}^{-1}$. In both cases, the fine-level correction is applied before the coarse-level correction.}
    
    \begin{enumerate}
    \item When $\mathbf{P}^{-1}_{\textrm{fine}} = \BI$, the preconditioner can be written as (according to Eq.~(19) in \cite{santo2018multi}):
    \begin{equation}
    \BM_{\textrm{V}1,l}^{-1} = \BI + \widetilde{\BU}^{(l)}_{r_{l}}\left((\widetilde{\BU}_{r_{l}}^{(l)})^{T}\WBA_{\bmu}\widetilde{\BU}^{(l)}_{r_{l}}\right)^{-1}(\widetilde{\BU}^{(l)}_{r_{l}})^{T}(\BI - \WBA_{\bmu}),
    \end{equation}
    where $\widetilde{\BU}_{r_{l}}^{(l)}$ is the reduced basis for the errors originating from the last fine-level iteration. 
    
    \item When using an SA preconditioner as the fine-level preconditioner, i.e., $\mathbf{P}_{\textrm{fine}}^{-1} = \BI + \BC^{-1}\BSigma_{s}$, the preconditioner becomes
    \begin{equation}
    \ning{\BM_{\textrm{V}2,l}^{-1}  = \BI+\BC^{-1}\BSigma_s +\widetilde{\BU}^{(l)}_{r_{l}}\left((\widetilde{\BU}_{r_{l}}^{(l)})^{T}\WBA_{\bmu}\widetilde{\BU}^{(l)}_{r_{l}}\right)^{-1}(\widetilde{\BU}^{(l)}_{r_{l}})^{T}(\BI - \WBA_{\bmu}(\BI+\BC^{-1}\BSigma_s)).}
    \end{equation}
    \end{enumerate}
These two direct extensions require the projection of the \ning{sweep-based scalar-flux operator}
$\WBA_{\bmu}=\BI-\BK_{\bmu}\BSigma_{s,\bmu}$, where
$\BK_{\bmu}=\sum_{j=1}^{N_{\Bupsilon}}
\omega_j(\BD_j+\BSigma_{t,\bmu})^{-1}$.
\ning{When the total cross section depends on the parameters,} the inverse transport operators generally introduce non-affine parameter dependence into $\WBA_{\bmu}$ \cite{behne2022minimally,peng2025flexible}. \ning{Consequently, assembling the reduced operator online requires full-order projections and $\sum_{l=1}^{N_w}r_l$ additional full-angle sweeps, where $N_w$ is the total aware level.} \ning{This online overhead can offset the convergence acceleration and even make the solver slower than DSA-preconditioned methods (see \cite[Sec. 4.2]{behne2022minimally}).}
\ning{Although hyper-reduction techniques such as EIM and DEIM~\cite{barrault2004empirical,chaturantabut2010nonlinear} can approximate such non-affine operators, evaluating sweep outputs at selected spatial points still demands upstream transport sweeps, limiting the computational savings. In contrast, our kinetic correction formulation allows the reduced matrix components to be precomputed offline under affine parameter dependence and naturally accommodates hyper-reduction under non-affine dependence. In both cases, online reduced-operator assembly of our approach avoids additional sweeps.}

\subsection{Practical implementation of ROM-based preconditioners}
In this section, we discuss the key algorithmic steps for efficiently implementing our method in practice.

\textbf{Offline precomputation:} The discrete operator for the ROM-based correction in the $l$-th iteration \pzc{with $1\leq l\leq N_w$} is defined as 
\begin{equation}
\BC_{\mathrm{ROM}, l}^{-1}  = (\sum_{j=1}^{N_{\Bupsilon}} \omega_{j}\BU^{(l)}_{r_{l},j})\left((\BU_{r_{l}}^{(l)})^{T}\BA_{\bmu}\BU^{(l)}_{r_{l}}\right)^{-1}(\sum_{j=1}^{N_{\Bupsilon}} \BU^{(l)}_{r_{l},j})^{T}.
\label{eq:rom-correction}
\end{equation}
Instead of assembling the matrix $\BC_{\mathrm{ROM}, l}^{-1}$ entirely during the online stage, we precompute and store the matrices 
\begin{equation}
    \BU_{r_l}^{\brho,(l)}=\sum_{j=1}^{N_{\Bupsilon}} \omega_{j}\BU^{(l)}_{r_{l},j} \quad\text{and}\quad
    \BU_{r_l}^{\textrm{iso},(l)}=\sum_{j=1}^{N_{\Bupsilon}} \BU^{(l)}_{r_{l},j}
\end{equation}
and the projected matrices $(\BU_{r_l}^{(l)})^T\BA_q\BU_{r_l}^{(l)}$ for the affine decomposition of $\BA_{\bmu}$ as discussed in Sec.~\ref{sec:rom}.

\textbf{Online application:} The crucial step to efficiently apply our SA preconditioner online is how to compute $\BC_{\mathrm{ROM}, l}^{-1}\BSigma_{s}\bq$ with $\bq\in\mathbb{R}^{N_\textrm{DOF}}$. Detailed steps are as follows. 
\begin{enumerate}
    \item Compute the matrix-vector product $\by=\BSigma_s\bq$ with $O(N_{\textrm{DOF}})$ cost.
    \item Compute \ning{$(\sum_{j=1}^{N_{\Bupsilon}}\BU_{r_l,j}^{(l)})^T\by$ as $(\BU_{r_l}^{\textrm{iso},(l)})^T\by$} with $O(rN_{\textrm{DOF}})$ cost.
    \item Solve the low-dimensional dense reduced-order problem $(\BU_{r_l}^{(l)})^T\BA_{\bmu}\BU_{r_l}^{(l)}\hat{\by}=\ning{(\BU_{r_l}^{\textrm{iso},(l)})^T\by}$ with $O(r^3)$ cost.
    \item Compute \ning{$(\sum_{j=1}^{N_{\Bupsilon}}\omega_j\BU_{r_l,j}^{(l)})\hat{\by}$ as $ \BU_{r_l}^{\rho,(l)}\hat{\by}$} with $O(rN_{\textrm{DOF}})$ cost.
\end{enumerate}
The total cost is $O(r(N_{\textrm{DOF}}+r^2))$. If one computes the matrix-matrix multiplication $\BC_{\mathrm{ROM}, l}^{-1}\BSigma_{s}$ first, the computational cost will rise to $O(N_{\textrm{DOF}}^2r)$.

\textbf{DSA:} After switching, we also need to apply the DSA preconditioner $\mathbf{M}^{-1}_{\mathrm{DSA}} = \BI +\BC_{\mathrm{DSA}}^{-1}\BSigma_{s}$, where \ning{$\BC_{\mathrm{DSA}}$} is the discrete diffusion operator. Specifically, we follow  \cite{Adams2002FastIM} to implement a fully or partially consistent discretization for the diffusion operator in DSA (see \cite[Appendix~B]{peng2024romsad} for details).


\section{Numerical results\label{sec:numerical}}
In this section, we present numerical results in both 1D and 2D settings to validate the efficiency and robustness of our trajectory-aware framework. We compare SI-TAR and FGMRES-TAR with SI-DSA, SI-ROMSAD, and GMRES-DSA, including their ROM-based initial-guess variants where applicable. For all simulations, we use linear upwind DG ($K=1$ in \eqref{eq:dg_space}) for spatial discretization. Unless otherwise specified, we apply a fully consistent DSA preconditioner derived \ning{by applying a $P_1$ moment closure to the discrete ideal correction equations and then eliminating the first-order moment variables using a Schur complement \cite{Adams2002FastIM}. The Robin boundary conditions are weakly imposed from the quadrature of transport numerical fluxes (see \cite[Appendix~B]{peng2024romsad} for details).}
The resulting diffusion system is solved using LU factorization in 1D slab geometry and the conjugate gradient method with an algebraic Multigrid (AMG) preconditioner in higher dimensions. The AMG solver is implemented using the
\texttt{iFEM} package \cite{chen2009integrated}. 

\ning{All methods use the following relative 2-norm residual as the stopping criterion:
\begin{equation}
\epsilon^{(l)}
=
\frac{\|\tilde{\bb}_{\bmu}
-\tilde{\BA}_{\bmu}\brho^{(l)}_{\bmu}\|_2}
{\|\tilde{\bb}_{\bmu}\|_2}
\label{eq:relative-residual}
\end{equation}
For SI, \eqref{eq:rho_residual} gives $\tilde{\bb}_{\bmu}
-\tilde{\BA}_{\bmu}\brho^{(l)}_{\bmu}=\brho^{(l+1,*)}_{\bmu}-\brho^{(l)}_{\bmu}$, which simplifies the stopping criterion to $\frac{\|\brho^{(l+1,*)}_{\bmu}-\brho^{(l)}_{\bmu}\|_2}
{\|\tilde{\bb}_{\bmu}\|_2}$,
where $\brho^{(l+1,*)}_{\bmu}$ is the intermediate density obtained from the $(l+1)$-th transport sweep before applying the SA correction.}

For those ROM-based methods, we define the  training set as $\mathscr{P}_{\mathrm{train}}$
with uniformly sampled parameters. 
We also randomly generate a test set $\mathscr{P}_{\mathrm{test}}$ \ning{using MATLAB's \texttt{rand} function with a fixed random seed initialized by \texttt{rng(10086)}}, ensuring that $\mathscr{P}_{\mathrm{train}} \cap \mathscr{P}_{\mathrm{test}} = \varnothing$. 

To analyze online performance on the test set, we calculate the average number of transport sweeps as:
\begin{equation}
    \bar{n}_{\textrm{sweep}} = \frac{\sum_{\bmu\in\mathscr{P}_{\textrm{test}}}\textrm{number of transport sweeps to reach  convergence for the parameter}\;\bmu }{\textrm{total number of parameters in }\mathscr{P}_{\textrm{test}}},
\end{equation}
and the average number of online iterations $\bar{n}_{\textrm{iter}}$ as referenced in Remark \ref{rem:gmres_sweep}. Then, we define the average $L_{\infty}$ residual: 
\begin{equation}
\bar{\mathcal{R}}_{\infty} = \frac{\sum_{\bmu\in\mathscr{P}_{\textrm{test}}}||(\BI-\mathbf{K}\BSigma_s)\brho_{\bmu}-{\wbb}_{\bmu}||_\infty}{\textrm{total number of parameters in }\mathscr{P}_{\textrm{test}}},
 \end{equation}
where $\brho_{\bmu}$ is the numerical solution from the underlying iterative solver. We also measure the average relative computational time $\bar{T}_{\textrm{rel}}$ with respect to SI-DSA. 

\ning{ For the 2D test cases, we further report two break-even numbers: $N_{\mathrm{BE}}$, which includes the full offline cost, and
$N_{\mathrm{BE}}^{\mathrm{reuse}}$, which assumes that high-fidelity training snapshots are already available:
\begin{equation}
N_{\mathrm{BE}} =
\left\lceil
\frac{T_{\mathrm{offline}}}
{\bar{T}_{\mathrm{SIDSA}}-\bar{T}_{\mathrm{online}}}
\right\rceil, \quad N_{\mathrm{BE}}^{\mathrm{reuse}} =
\left\lceil
\frac{T_{\mathrm{offline}}-T_{\mathrm{snapshot}}}
{\bar{T}_{\mathrm{SIDSA}}-\bar{T}_{\mathrm{online}}}
\right\rceil
\end{equation}
where $T_{\mathrm{offline}}$ denotes the total offline cost of the corresponding method, including the snapshot generation cost
$T_{\mathrm{snapshot}}$. The quantities
$\bar{T}_{\mathrm{SIDSA}}$ and $\bar{T}_{\mathrm{online}}$ denote the average solution times per parameter for SI-DSA and the corresponding online method, respectively. For each method, both averages are computed from timings recorded in that method’s own experimental run.
}

\ning{All methods were implemented in MATLAB R2026a. All numerical experiments were conducted on a 2024 MacBook Air equipped with an Apple M3 chip and 16 GB of unified memory.}

To improve readability, we summarize the abbreviations used throughout this section in Table~\ref{tab:abbreviation}.

\begin{table}[!ht]
\centering
\caption{Abbreviations used in Sec.~\ref{sec:numerical}.}
\label{tab:abbreviation}
\footnotesize

\begin{subtable}[t]{\textwidth}
\centering
\renewcommand{\arraystretch}{1.2}
\begin{tabular}{cl}
\toprule
\multirow{2}{*}{SI-ROMSAD-$\mathfrak{w},l$}
& ROMSAD preconditioner \cite{peng2024romsad,peng2025flexible} using window size $\mathfrak{w}$ to build the ROM \\
& and switching to DSA after the $l$-th iteration \\
\midrule

SI-TAR(-$N_w$)
& SI with trajectory-aware ROM-enhanced SA and 0 initial guess (using $N_w$ aware levels) \\
\midrule

\multirow{2}{*}{SI-TAR-IG(-$N_w$)}
& SI with trajectory-aware ROM-enhanced SA \\
& and a ROM-based initial guess (using $N_w$ aware levels) \\
\midrule

GMRES-DSA
& GMRES with right DSA preconditioner and 0 initial guess \\
\midrule

GMRES-DSA-IG
& GMRES with right DSA preconditioner and a ROM-based initial guess \\
\midrule

\multirow{2}{*}{\ning{FGMRES-TAR-IG}(-$N_w$)}
& FGMRES with trajectory-aware ROM-enhanced SA \\
& and a ROM-based initial guess (using $N_w$ aware levels) \\
\bottomrule
\end{tabular}
\end{subtable}
\end{table}

\begin{rem}
In the FGMRES-based method (Alg.~\ref{alg:fgmres}), computing the initial residual $\br^{(0)}=\wbb-\WBA\brho^{(0)}$ requires one additional transport sweep to evaluate $\wbb$ for a zero initial guess, and a further sweep to evaluate $\WBA\brho^{(0)}$ for a nonzero initial guess. \ning{Similarly, SI-based methods require one additional sweep to evaluate the denominator of the normalized residual stopping criterion. All these sweeps are included in the reported totals.}
\label{rem:gmres_sweep}
\end{rem}

\subsection{Two-material problem in 1D slab geometry\label{sec:two_material}}
We start from a 1D slab geometry configuration: 
\begin{equation}
    \xi\partial_x f(x,\xi)+\sigma_t(x) f(x,\xi)=\sigma_s(x)\rho(x)+G(x),
\end{equation}
where $\xi=\cos(\theta)\in[-1,1]$.
The goal is to investigate how hyperparameters, such as the POD truncation tolerance $\epsilon_{\mathrm{POD}}$ and the total aware level $N_{w}$, affect the performance of the proposed methods. 
Consider a parametric two-material problem as follows:
\begin{subequations}
\label{eq:two_material}
    \begin{align}
G(x)= 0 ,\quad \sigma_{a} = \left\{
\begin{array}{ll}
\mu_a, & 0\leq x\leq 1, \\
0, & 1<x\leq11,
\end{array}
\right. \quad 
\ning{\sigma_{s}} = \left\{
\begin{array}{ll}
0, & 0\leq x\leq 1, \\
\mu_{s}, & 1<x\leq11,
\end{array}
\right. \quad 
x \in[0,11],
\end{align}
\begin{align}
    f(0,\xi) = 5 \;\;\textrm{with}\;\; \xi>0, \quad f(11,\xi) = 0 \;\;\textrm{with}\;\; \xi<0.
\end{align}
\end{subequations}
The computational domain is partitioned into two distinct regions: a pure absorption zone $(x\in[0,1])$ with \ning{an} absorption cross section $\mu_{a}\in[0.5,1.5]$, and a pure scattering zone $(x\in(1,11])$ with a scattering cross section $\mu_{s}\in[10, 50]$. For spatial discretization, we employ a non-uniform mesh using a refined grid $\Delta x = \Delta x_1 = \frac{1}{100}$ in the absorption region $[0,1]$ and a coarse grid $\Delta x = \Delta x_2 = \frac{1}{10}$ in the scattering region $(1,11]$. We use a $16$-point Gauss-Legendre quadrature over $[-1,1]$ for angular discretization.

\begin{figure}[!ht]
    \centering
    \begin{subfigure}[b]{0.4\textwidth}
        \centering
        \begin{minipage}{6cm}
        \includegraphics[scale=0.5]{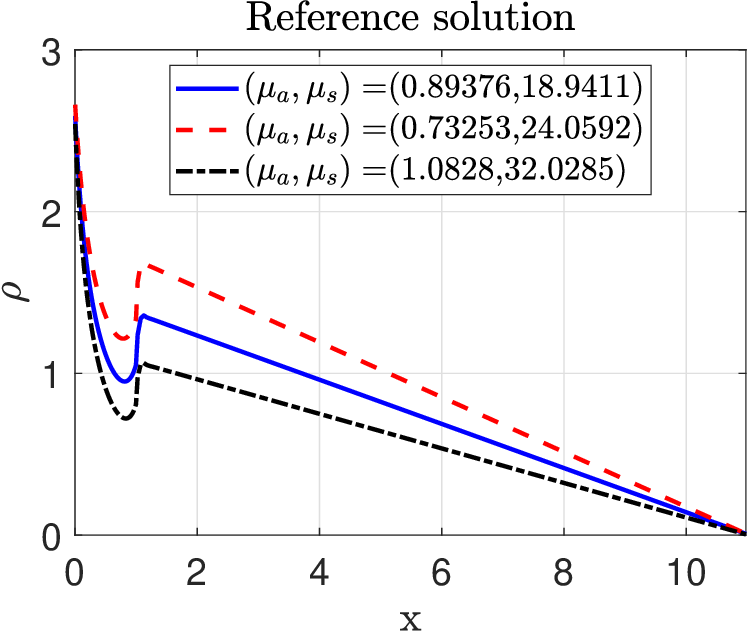}
        \end{minipage}
        \label{fig:two-material-sol}
    \end{subfigure}
    \begin{subfigure}[b]{0.45\textwidth}
        \centering
        \begin{minipage}{6cm}
        \includegraphics[scale=0.5]{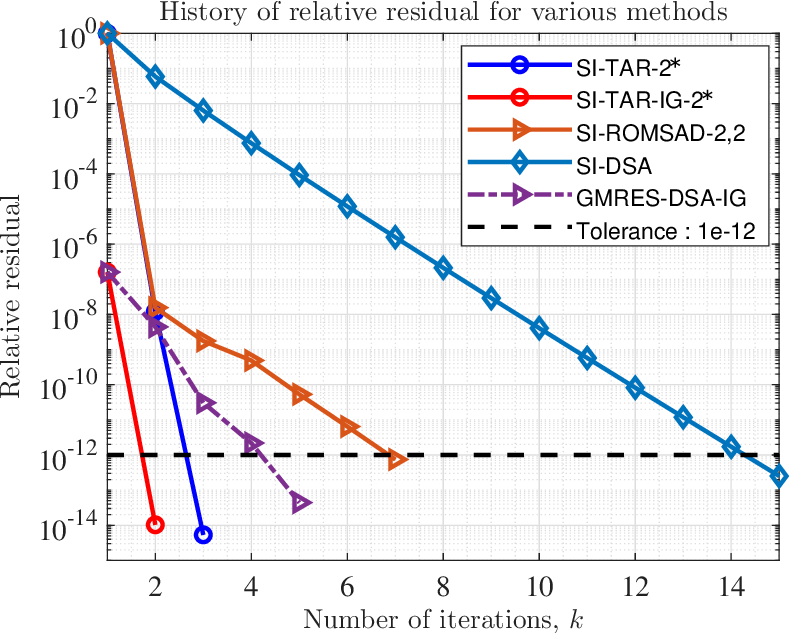}
        \end{minipage}
        \label{fig:two-material-his}
    \end{subfigure}
    \caption{Results for the 1D two-material problem in Sec.~\ref{sec:two_material}. Left: reference solution. Right: relative residual histories for the parameter pair $(\mu_{a},\mu_s)=(0.73253,24.0592)$. 
    } 
    \label{fig:two_material}
\end{figure}

The training set for this problem is constructed as follows:
\begin{equation}
    \mathscr{P}_{\mathrm{train}} = \{\;(\mu_{a,m}, \mu_{s,n}) \;:\; \mu_{a,m} = 0.5+0.1m, \;\mu_{s,n} = 10+n,\;\; 0\leq m\leq 10,\; 0\leq n\leq 40\;\},
\end{equation}
where $m$ and $n$ are integer indices. To evaluate the performance of the proposed methods, we randomly select 20 parameter pairs $(\mu_a, \mu_s)$ from the parameter space $[0.5,1.5]\times[10,50]$ as test cases.  The reference solutions corresponding to three different parameter pairs $(\mu_a, \mu_s)$ are visualized on the left-hand side of \ning{Fig.~\ref{fig:two_material}}.

Table~\ref{tab:two_material} summarizes the computational results for \ning{two} sets of hyperparameters,  $\epsilon_{\mathrm{POD}}$ and total aware level $N_{w}$, where the tolerance for iteration convergence is $\ning{\epsilon_{\textrm{SI}}}=\epsilon_{\textrm{GMRES}}=10^{-12}$. We have the following key observations.

\textbf{Influence of POD tolerance $\epsilon_{\mathrm{POD}}$:} As $\epsilon_{\mathrm{POD}}$ decreases, trajectory-aware ROM-based methods converge with fewer transport sweeps, since the underlying ROMs become more accurate. 

\textbf{Influence of total aware level $N_{w}$:} 
\ning{For both SI- and FGMRES-based TAR methods, higher total aware level $N_{w}$ can reduce the number of transport sweeps. With $\epsilon_{\mathrm{POD}}:  = 10^{-5}$, increasing $N_w$ from 1 to 2 results in an average 1.45-iteration reduction for FGMRES-TAR-IG. This observation implies that the trajectory-aware framework effectively approximates the ideal correction in the second iteration. Moreover, when $\epsilon_{\mathrm{POD}}:  = 10^{-7}$, a single aware level is sufficient for FGMRES-TAR-IG to achieve convergence within one iteration.}

\textbf{Online performance compared to DSA and ROMSAD:} Our main observations are as follows.
\begin{enumerate}
\item
To visually illustrate the benefit of the trajectory-aware framework, we present the residual history for  $(\mu_{a},\mu_s)=(0.73253,24.0592)$ in Fig.~\ref{fig:two_material}. Using two aware levels, the trajectory-aware ROM mitigates the efficiency loss observed in ROMSAD, enabling \ning{SI-TAR-IG-2 to converge in only two iterations.} 
\item
With a zero initial guess and $\epsilon_{\mathrm{POD}}=10^{-5}$, SI-TAR-2 reduces transport sweeps by factors of approximately \ning{2.18} and \ning{1.61} relative to SI-DSA and SI-ROMSAD-2,2, respectively. At $\epsilon_{\mathrm{POD}}=10^{-7}$, SI-TAR-1 performs comparably to SI-ROMSAD-2,2, \ning{while introducing a second aware level reduces the number of sweeps by factors of approximately $3.88$ and $2.20$ compared to SI-DSA and SI-ROMSAD-2,2, respectively.}
\item When only one aware level is applied, using a ROM-based initial guess for SI significantly accelerates the convergence.
\item
When using a ROM-based initial guess and $\epsilon_{\textrm{POD}}=10^{-5}$, \ning{FGMRES-TAR-IG}-$1$ and \ning{FGMRES-TAR-IG}-$2$ reduce the number of transport sweeps for convergence by factors of approximately \ning{$1.25$} and \ning{$1.70$} compared to \ning{GMRES-DSA-IG}. When lowering $\epsilon_{\textrm{POD}}$ to $10^{-7}$, \ning{FGMRES-TAR-IG} solvers with both $N_w=1$ and $N_w=2$ converge after one iteration and reduce the required number of transport sweeps by a factor of approximately \ning{$1.87$}.
\end{enumerate}
The primary objective of this 1D slab geometry test is to examine the impact of hyperparameters. In fact, as demonstrated in later tests, more significant online savings can be achieved by our method in higher dimensions.

\begin{table}[!ht]
\small
\centering
\ning{
\caption{Results for the 1D two-material problem (Sec.~\ref{sec:two_material}). $r_{\textrm{IG}}$: the dimension of the ROM for the initial guess. \label{tab:two_material}} 
\vspace{-10pt}
\begin{subtable}[t]{1.0\textwidth}
\centering
\caption{Results for various SI methods. Denote $r_{c}$ as the dimension of the ROM for the correction equation. \ning{When $\epsilon_{\mathrm{POD}}=10^{-7}$, for SI-TAR-1, $r_c = 25$, for SI-TAR-IG-1, $r_c = 88$, for SI-TAR-2, $r_{c,1} =25, r_{c,2}=85$, for SI-TAR-IG-2, $r_{c,1} = 88 , r_{c,2} = 451$.}}
\vspace{8pt}
\begin{tabular}{cccccccc}
\toprule
 & SI-DSA & SI-ROMSAD-2,2 & SI-TAR-1 & SI-TAR-IG-1 & SI-TAR-2 & SI-TAR-IG-2 \\
\midrule
\multicolumn{7}{c}{(1) $\epsilon_{\mathrm{POD}} = 10^{-5}$, $r_{\textrm{IG}} =13$ 
}\\
\midrule
$\bar{n}_{\textrm{sweep}}$ & $15.50$ & $11.40$ & $10.70$ & $5.80$ & $7.10$ & $4.00$\\
$\bar{\mathcal{R}}_{\infty}$& $5.60\times10^{-13}$& $8.31\times10^{-13}$& $1.01\times10^{-12}$& $1.31\times10^{-12}$& $7.33\times 10^{-13}$ & $1.05\times 10^{-13}$ \\
\midrule
\multicolumn{7}{c}{(2) $\epsilon_{\mathrm{POD}}  = 10^{-7}$, $r_{\textrm{IG}} =25$}\\
\midrule
$\bar{n}_{\textrm{sweep}}$& $15.50$ & $8.80$ & $8.85$ & $3.00$ &  $4.00$ & $3.00$  \\
$\bar{\mathcal{R}}_{\infty}$& $5.60\times10^{-13}$& $9.31\times10^{-13}$& $9.40\times10^{-13}$ & $8.40\times10^{-14}$ & $1.50\times 10^{-13}$ &  $8.40\times 10^{-14}$\\
\bottomrule
\end{tabular}
\end{subtable}
\\[1em]
\begin{subtable}{1.0\textwidth}
\small
\centering
\caption{Results for various GMRES methods. Denote $r_{c}$ as the dimension of the ROM for the correction equation. When $\epsilon_{\mathrm{POD}}=10^{-7}$, \ning{for FGMRES-TAR-IG-1, $r_c = 89 $, for FGMRES-TAR-IG-2, $r_{c,1} =89, r_{c,2} = 441$. }}
\vspace{8pt}
\begin{tabular}{cccccccccc}
\toprule
  &  GMRES-DSA & GMRES-DSA-IG & \ning{FGMRES-TAR-IG}-1 & \ning{FGMRES-TAR-IG}-2 \\
\midrule
\multicolumn{5}{c}{(1) $\epsilon_{\mathrm{POD}} = 10^{-5}$, $r_{\textrm{IG}} =13$}\\
\midrule
$\bar{n}_{\textrm{iter}}$ & $8.00$ & $4.80$ & $3.45$ & $2.00$ \\
$\bar{n}_{\textrm{sweep}}$ & $9.00$ & $6.80$ & $5.45$ & $4.00$\\
$\bar{\mathcal{R}}_{\infty}$& $1.22\times10^{-12}$& $1.91\times10^{-12}$ & $1.02\times 10^{-12}$ & $1.73\times10^{-13}$\\
\midrule
\multicolumn{5}{c}{(2) $\epsilon_{\mathrm{POD}}  = 10^{-7}$, $r_{\textrm{IG}} =25$}\\
\midrule
$\bar{n}_{\textrm{iter}}$ & $8.00$ &$3.60$ & $1.00$ & $1.00$ \\
$\bar{n}_{\textrm{sweep}}$& $9.00$ &$5.60$ & $3.00$ & $3.00$ \\
$\bar{\mathcal{R}}_{\infty}$&  $1.22\times10^{-12}$ &$1.55\times10^{-12}$ & $7.81\times10^{-14}$ & $7.81\times10^{-14}$\\
\bottomrule
\end{tabular}
\end{subtable}
}
\end{table}


\subsection{Variable scattering problem\label{sec:scattering}}
\begin{figure}[!ht]
    \centering
    \begin{subfigure}[b]{0.4\textwidth}
        \centering
        \begin{minipage}{6cm}
        \includegraphics[scale=0.5]{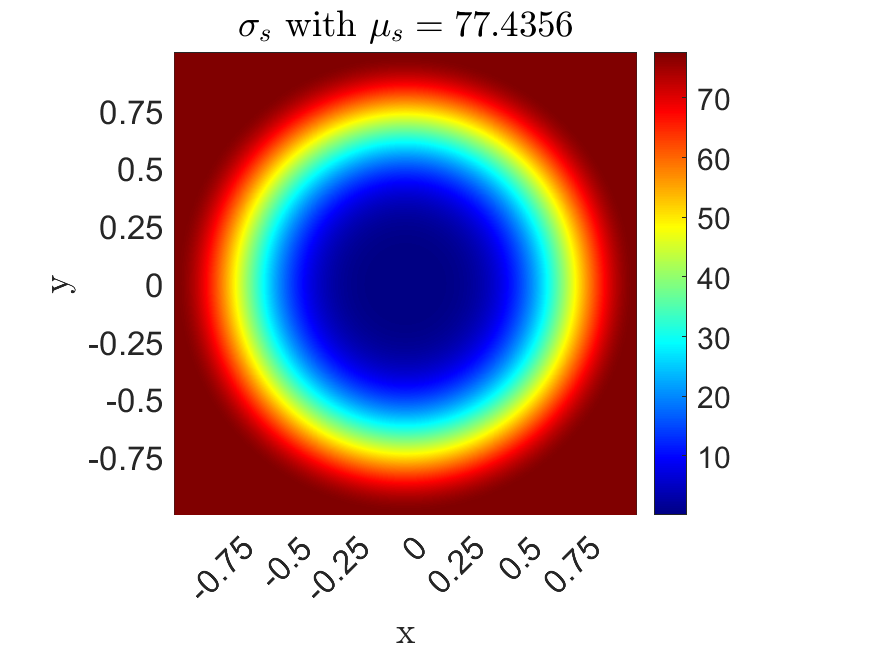}
        \end{minipage}
        \label{fig:scattering-var}
    \end{subfigure}
     \hspace{-0.04\textwidth}
    \begin{subfigure}[b]{0.45\textwidth}
        \centering
        \begin{minipage}{6cm}
        \includegraphics[scale=0.5]{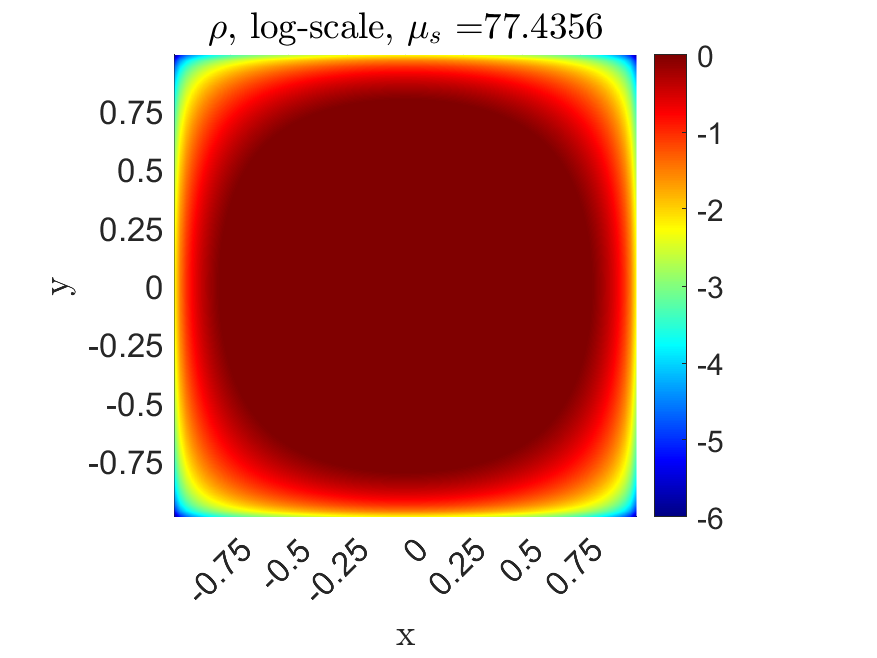}
        \end{minipage}
        \label{fig:scattering-ref}
    \end{subfigure}
    \caption{ Left: the setup of the scattering cross section for the  variable scattering problem in Sec.~\ref{sec:scattering}. Right: a reference solution.} 
    \label{fig:scattering_variable}
\end{figure}

We consider a parametric problem with spatially variable scattering on the computational domain $\mathcal{D} = [-1,1]^{2}$, where the absorption cross section $\sigma_a$ is zero, and the scattering cross section is defined as
\begin{equation}
   \sigma_s(x,y) = \left\{
\begin{array}{ll}
\mu_sr^{4}(2-r^2)^2 + 0.1, \quad & r = \sqrt{x^2+y^2}\leq 1, \\
\mu_s+0.1, & \textrm{otherwise},
\end{array}
\right.
\quad \mu_{s}\in[49.9,99.9].
\end{equation}
This scattering cross section varies smoothly from 0.1 at the domain center to $\mu_{s} + 0.1 \in[50,100]$ at the boundary. This variation indicates a gradual transition from a transport-dominated regime to a scattering-dominated regime. The parameter $\mu_{s}$ controls the rate of this transition. The problem is solved with zero-inflow boundary conditions and a Gaussian source term given by $G(x,y) = \frac{10}{\pi}\exp(-100(x^2+y^2))$. In Fig.~\ref{fig:scattering_variable}, we present the scattering cross section $\sigma_{s}(x,y)$ with \ning{$\mu_{s} = 77.4356$} and the corresponding reference solution. We use an $80\times80$ uniform mesh for spatial discretization and the CL(30,6) quadrature rule for angular discretization. We set the convergence tolerances for SI and FGMRES to $\ning{\epsilon_{\textrm{SI}}}=\epsilon_{\mathrm{GMRES}}=10^{-11}$.

\ning{This test establishes how resolving the offline-online trajectory mismatch improves correction accuracy and overall online  efficiency. We also assess how fast the resulting online speedups amortize the additional offline construction cost. }
We construct ROMs using a training dataset $\mathscr{P}_{\mathrm{train}}$ comprising 50 uniformly distributed samples. For testing purposes, we randomly sample 10 values of $\mu_{s}$ from $[49.9, 99.9]$.

\textbf{Online performance:} We consider two configurations of the POD truncation tolerance and the total aware level: $(\epsilon_{\mathrm{POD}},N_w)=(10^{-5},2)$ and $(10^{-7},1)$. During the online phase, we evaluate the performance of the proposed methods using 10 randomly selected test samples. As shown in \ning{Table~}\ref{tab:changing_scattering_online}, all trajectory-aware methods demonstrate significant acceleration over both SI-DSA and DSA-preconditioned GMRES:
\begin{enumerate}
    \item When $\epsilon_{\mathrm{POD}}=10^{-5}$ and $N_w=2$, SI-TAR-IG reduces the number of transport sweeps \ning{by factors of approximately $4.25$ and $2.80$} compared to SI-DSA and \ning{SI-ROMSAD-IG-2,2}, respectively, achieving corresponding speedups of \ning{$7.44$ and $4.50$}. When $\epsilon_{\mathrm{POD}}=10^{-7}$ and $N_w=1$, SI-TAR-IG converges in only $3$ transport sweeps: \ning{$2$ SI online sweeps, and $1$ RHS sweep to compute the normalization factor $\|\tilde{\bb}\|_2$}. This yields sweep reduction factors of \ning{$5.67$ and $2.43$, with corresponding computational speedups of $11.20$ and $3.84$}, relative to SI-DSA and SI-ROMSAD-IG-2,2, respectively.
 
    \item At $\epsilon_{\mathrm{POD}}=10^{-7}$,  \ning{FGMRES-TAR-IG} requires an average of only three transport sweeps, reducing the number of transport sweeps by factors of approximately 3.67 and 1.67 compared to GMRES-DSA and GMRES-DSA-IG, respectively.

 \item \ning{The two TAR-IG variants achieve greater wall-clock speedups than their sweep reductions alone suggest because they converge within the ROM-preconditioned iterations without switching to DSA. In contrast, for $\epsilon_{\mathrm{POD}}=10^{-5}$ and $N_w=2$, DSA solves add $40.73\%$, $28.89\%$, $42.51\%$, and $33.30\%$ to the sweep time in SI-DSA, SI-ROMSAD-IG-2,2, GMRES-DSA, and GMRES-DSA-IG, respectively.}

\end{enumerate}
Among the methods tested, SI-TAR-IG and FGMRES-TAR-IG achieve comparable online performance and outperform the others, which emphasizes the importance of accounting for preconditioner dependence in reduced-order modeling.

\begin{table}[!ht]
\centering
\ning{
\caption{Test results for the 2D variable scattering problem (Sec.~\ref{sec:scattering}). Denote $N_w$ as the total aware level, $r_{\mathrm{IG}}$ as the dimension of the ROM for the initial guess, and $r_{c,l}$ as the dimension of the ROM at the $l$-th aware level for the correction equation. \ning{When $\epsilon_{\mathrm{POD}}=10^{-7}$ and $N_w=1$, $r_{c,1}=15$ for SI-TAR-IG and $r_{c,1}=20$ for FGMRES-TAR-IG.}}
\vspace{8pt}
\label{tab:changing_scattering_online}
\begin{subtable}[t]{\textwidth}
\centering
\small
\setlength{\tabcolsep}{4pt}
\begin{adjustbox}{max width=\linewidth,center}
\begin{tabular}{@{}lcccccc@{}}
\toprule
 & SI-DSA
 & \ning{SI-ROMSAD-IG-2,2}
 & SI-TAR-IG
 & GMRES-DSA
 & GMRES-DSA-IG
 & \ning{FGMRES-TAR-IG} \\
\midrule
\multicolumn{7}{c}{(1) $\epsilon_{\mathrm{POD}}=10^{-5}$,
$r_{\mathrm{IG}}=4$, $N_w=2$}\\
\midrule
$\bar{n}_{\textrm{sweep}}$
& 17 & 11.2 & 4 & 11 & 7 & 4 \\
$\bar{\mathcal{R}}_{\infty}$
& $2.04\times10^{-12}$
& $5.18\times10^{-13}$
& $2.72\times10^{-13}$
& $3.21\times10^{-13}$
& $1.03\times10^{-12}$
& $3.21\times10^{-13}$ \\
$\bar{T}_{\textrm{rel}}$
& 100\% & 60.43\% & 13.44\%
& 62.57\% & 37.23\% & 13.08\% \\
$N_{\mathrm{BE}}\;(N_{\mathrm{BE}}^{\mathrm{reuse}})$
& -- & 206\;(49) & 80\;(8)
& -- & -- & 87\;(16) \\
\midrule
\multicolumn{7}{c}{(2) $\epsilon_{\mathrm{POD}}=10^{-7}$,
$r_{\mathrm{IG}}=6$, $N_w=1$}\\
\midrule
$\bar{n}_{\textrm{sweep}}$
& 17 & 7.3 & 3 & 11 & 5 & 3 \\
$\bar{\mathcal{R}}_{\infty}$
& $2.04\times10^{-12}$
& $4.20\times10^{-13}$
& $1.82\times10^{-13}$
& $3.21\times10^{-13}$
& $9.02\times10^{-13}$
& $5.15\times10^{-13}$ \\
$\bar{T}_{\textrm{rel}}$
& 100\% & 34.32\% & 8.93\%
& 65.15\% & 22.91\% & 8.50\% \\
$N_{\mathrm{BE}}\;(N_{\mathrm{BE}}^{\mathrm{reuse}})$
& -- & 119\;(28) & 70\;(4)
& -- & -- & 74\;(9) \\
\bottomrule
\end{tabular}
\end{adjustbox}
\end{subtable}
}
\end{table}

\textbf{\ning{Offline-online mismatch mitigation through trajectory-aware ROMs.}} \pzc{
To explain the significant online speedup gained by the trajectory-aware approach,
we isolate the effect of the offline-online trajectory mismatch by comparing the projection errors of TAR against ROMSAD and a
\ning{fixed-trajectory multispace ROM (FTMS)}, which constructs iteration-dependent reduced spaces \ning{from offline trajectories generated with fixed SI-DSA.} Hence, FTMS does not account for how earlier ROM corrections alter the subsequent online trajectory. \ning{Table~\ref{tab:projection_errors} reports the mean relative $\ell_2$ projection error of the online ideal correction $\delta\bff_{\mu,\mathrm{ideal}}^{(l)}$ onto the reduced space spanned by the orthonormal basis $\BU_r^{(l)}$, averaged over the test samples for each method:}}
\begin{equation}
    e_{\mathrm{proj}}^{(l)}(\mu)
    = {\left\|\left(\BI - \BU_r^{(l)}(\BU_r^{(l)})^T\right)\delta\bff_{\mu,\mathrm{ideal}}^{(l)}\right\|_2}/{\left\|\delta\bff_{\mu,\mathrm{ideal}}^{(l)}\right\|_2}.
\end{equation}

\ning{At the first aware level ($l=1$), most methods achieve comparable projection accuracy because their offline and
online initial-guess settings are consistent.
The only exception is SI-ROMSAD-IG-2,2, whose snapshots are collected
from zero-initial-guess solves although its online solves use a ROM-based initial guess}.

\ning{At the second aware level ($l=2$), however, the first
ROM correction changes the solution trajectory, and TAR achieves higher projection accuracy by accounting for
this change during offline construction.
Specifically, TAR reduces the mean projection error relative to FTMS
by factors of approximately $2.8\times10^4$ and $1.1\times10^2$ for SI, and $2.9\times10^4$ and $10$ for FGMRES, under zero and ROM-based initial guesses, respectively. SI-TAR and SI-TAR-IG similarly outperform SI-ROMSAD-2,2 and SI-ROMSAD-IG-2,2 by factors of approximately $1.3\times10^3$ and $29$,
respectively}.

\begin{table}[!ht]
\centering
\caption{Mean relative $\ell_2$ projection errors of the ideal angular-flux corrections for the 2D variable-scattering problem ($\varepsilon_{\mathrm{POD}}=10^{-5},\;N_w=2$) with zero and ROM-based initial guesses. The latter correspond to the -IG variants. }
\label{tab:projection_errors}
\ning{
\begin{tabular}{lcccc}
\toprule
& \multicolumn{2}{c}{Zero initial guess}
& \multicolumn{2}{c}{ROM initial guess} \\
\cmidrule(lr){2-3}
\cmidrule(lr){4-5}
Method
& $l=1$ & $l=2$
& $l=1$ & $l=2$ \\
\midrule
SI-TAR
& $9.47\times10^{-6}$
& $1.87\times10^{-5}$
& $5.24\times10^{-6}$
& $4.29\times10^{-3}$ \\
SI-FTMS
& $9.47\times10^{-6}$
& $5.24\times10^{-1}$
& $5.24\times10^{-6}$
& $4.80\times10^{-1}$ \\
SI-ROMSAD-2,2
& $8.12\times10^{-6}$
& $2.41\times10^{-2}$
& $1.06\times10^{-1}$
& $1.26\times10^{-1}$ \\
\midrule
FGMRES-TAR
& $9.44\times10^{-6}$
& $3.60\times10^{-6}$
& $3.74\times10^{-6}$
& $2.59\times10^{-3}$ \\
FGMRES-FTMS
& $9.44\times10^{-6}$
& $1.04\times10^{-1}$
& $3.74\times10^{-6}$
& $2.59\times10^{-2}$ \\
\bottomrule
\end{tabular}}
\end{table}

\textbf{Offline efficiency \pzc{and break-even analysis}:} The offline cost of the trajectory-aware framework consists of three main components: (i)  SI-DSA solves to generate full-order solutions for 50 training parameters;  (ii) one additional full-angle transport sweep per parameter at each aware level; (iii) construction of reduced bases and projected operators.
Table~\ref{tab:offline_time_scattering} reports the relative additional offline computational costs besides snapshot generation for $\epsilon_{\mathrm{POD}}=10^{-7}$ and $N_w=1$, normalized by a single SI-DSA linear solve. The transport-sweep cost $\bar{T}_{\mathrm{sweep}}$ and Arnoldi cost for tracking the FGMRES trajectory $\bar{T}_{\mathrm{Arnoldi}}$ are averaged over training samples, \ning{whereas the basis construction cost $T_{\mathrm{basis}}$ and operator construction cost $T_{\mathrm{operator}}$ are reported as total costs, since they are constructed collectively from the entire training set rather than separately for each parameter.} Our main observations are as follows.

\ning{As shown in Table~\ref{tab:offline_time_scattering}, the dominant offline cost comes from the $50$ full-order SI-DSA solves to generate training data, while the additional ROM-related computations incur only modest overhead relative to these training solves. Specifically, building the ROM-based initial guess requires only $3.30\%$ ($T_{\mathrm{basis}}$) and $0.28\%$ ($T_{\mathrm{operator}}$) of a single SI-DSA solve. For the correction ROMs, basis assembly, operator projection, and snapshot sweeps require $3.94\%$, $0.65\%$, and $4.88\%$ under SI, and $3.46\%$, $0.84\%$, and $4.10\%$ under FGMRES (plus $4.14\%$ for the Arnoldi process). These additional sweeps and Arnoldi updates beyond
ROMSAD~\cite{peng2024romsad,peng2025flexible}
add negligible cost relative to the solves generating training data.}

\ning{Given the modest offline overhead and substantial online acceleration, the trajectory-aware framework achieves low break-even counts ($N_{\mathrm{BE}}$ in Table~\ref{tab:changing_scattering_online}). For $\epsilon_{\mathrm{POD}}=10^{-5}$ and $10^{-7}$, FGMRES-TAR-IG breaks even after $87$ and $74$ online evaluations, respectively. When high-fidelity training snapshots are reused ($N_{\mathrm{BE}}^{\mathrm{reuse}}$), these counts decrease to $16$ and $9$ evaluations for FGMRES-TAR-IG and to $8$ and $4$ for SI-TAR-IG. Compared to SI-ROMSAD-IG-2,2, SI-TAR-IG reduces the break-even count by factors of $2.58$ and $1.70$ under full offline costs (and from $49$ and $28$ down to $8$ and $4$ when reusing snapshots).}

\begin{table}[!ht]
\centering
\ning{
\begin{threeparttable}
\captionsetup{
    labelformat=simple,
    labelsep=colon
}
\caption{The average relative offline computational cost for the 2D variable scattering problem (Sec.~\ref{sec:scattering}), where $\epsilon_{\mathrm{POD}} = 10^{\ning{-7}}, N_w=1$. \ning{A dash ``--'' indicates that the entry is not applicable.}}
\label{tab:offline_time_scattering}
\begin{tabular}{lcccc}
\toprule
 & $T_{\mathrm{basis}}$ & $T_{\mathrm{operator}}$ & $\bar{T}_{\mathrm{sweep}}$ & $\bar{T}_{\mathrm{Arnoldi}}$\\
\midrule
Initial guess
& 3.30\% & 0.28\% & -- & -- \\
Correction equation (SI-TAR-IG)
& 3.94\% & 0.65\% & 4.88\% & -- \\
Correction equation (FGMRES-TAR-IG)
& 3.46\% & 0.84\% & 4.10\% & 4.14\% \\
\bottomrule
\end{tabular}
\end{threeparttable}}
\end{table}

\subsection{Pin-cell problem\label{sec:pincell}}
\begin{figure}[!ht]
    \centering
    \begin{subfigure}[b]{0.4\textwidth}
        \centering
        \begin{minipage}{6cm}
        \includegraphics[scale=0.5]{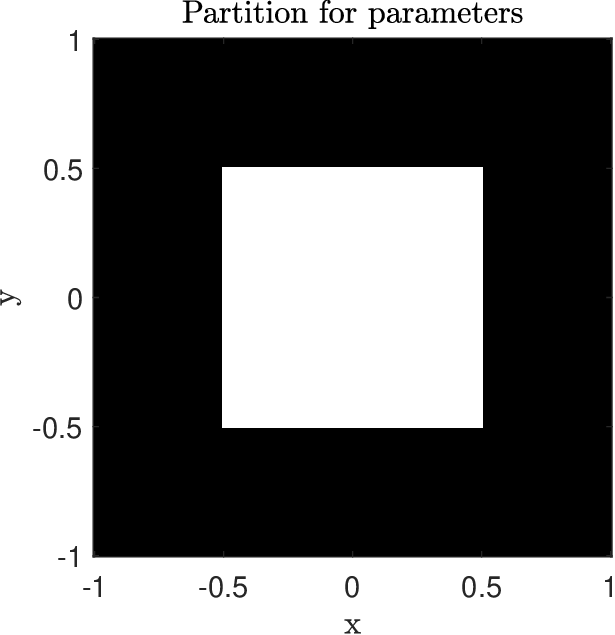}
        \end{minipage}
        \label{fig:pincell-var}
    \end{subfigure}
    \hspace{-0.04\textwidth}
    \begin{subfigure}[b]{0.4\textwidth}
        \centering
        \begin{minipage}{6cm}
        \includegraphics[scale=0.5]{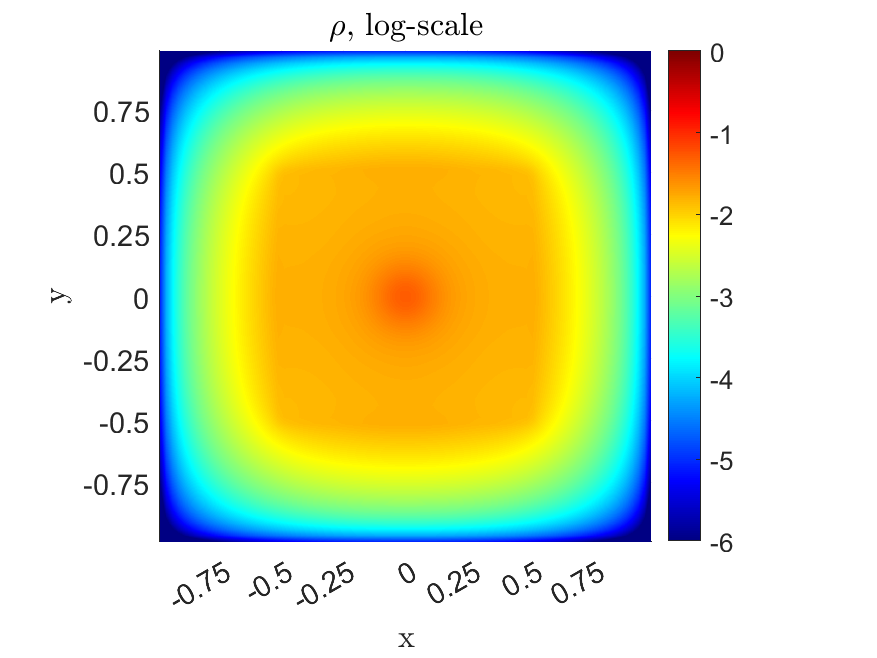}
        \end{minipage}
        \label{fig:pincell-ref}
    \end{subfigure}
    \caption{ Left: the setup for the pin-cell problem in Sec.~\ref{sec:pincell}, where the outer black region represents a pure scattering medium with $\sigma_s = 100$, while the inner white region has parametric scattering and absorption cross sections. Right: reference solution for $(\mu_a,\mu_s) = (0.1506,0.2978)$.} 
    \label{fig:pincell_variable}
\end{figure}

We consider a parametric pin-cell problem with zero inflow boundary conditions on the computational domain $\mathcal{D} = [-1,1]^{2}$, where the source term is given by $G(x,y) = \exp(-100(x^2+y^2))$. The left panel of Fig.~\ref{fig:pincell_variable} shows the material configuration. The outer black region represents a purely scattering medium with $\sigma_s = 100$, while the inner white region has parametric scattering and absorption cross sections defined by:
\begin{equation}
(\sigma_a, \sigma_s) = 
\left\{
\begin{array}{ll} 
    (\mu_a, \mu_s)\in[0.05,0.5]^2, &\quad \text { if }|x| \leq 0.5 \text{ and }|y| \leq 0.5, \\
    (0, 100), &\quad \text { otherwise. }\end{array} 
\right. 
\end{equation}
This setup introduces a dramatic scattering transition between adjacent regions, resulting in challenging multiscale effects. We use an $80\times80$ uniform spatial grid with CL(30,6) angular quadrature, and employ partially consistent DSA for acceleration. The convergence tolerance is fixed at $\ning{\epsilon_{\textrm{SI}}}=\epsilon_{\mathrm{GMRES}}=10^{-11}$.

\ning{This test evaluates trajectory-aware acceleration in a challenging
multiscale problem, where SI-DSA fails to reach the prescribed tolerance within 50 transport sweeps. We examine how the ROM configuration and initial guess affect
convergence and efficiency.}

The training set for ROM construction consists of 25 uniformly distributed samples:
\begin{equation}
    \mathscr{P}_{\text {train }}=\left\{
    \left(\mu_a, \mu_s\right)=\left(\ning{0.05}+i \Delta \mu_a, \ning{0.05}+j \Delta \mu_s\right),\;\Delta \mu_a=\Delta \mu_s=\ning{0.1125,}\;\ning{i, j=0,\ldots,4}\right\}.
\end{equation}
For online testing, we randomly select an additional 10 parameter pairs $(\mu_a, \mu_s)\in[0.05,0.5]^2$. The reference solution at $(\mu_a,\mu_s) = (0.1506,0.2978)$ is shown in Fig.~\ref{fig:pincell_variable}.

\begin{table}[!ht]
\centering
\ning{\caption{\ning{Test results for the 2D pin-cell problem (Sec.~\ref{sec:pincell}). Denote $N_w$ as the total aware level, $r_{\mathrm{IG}}$ as the dimension of the ROM for the initial guess, and $r_{c,l}$ as the dimension of the ROM at the $l$-th aware level for the correction equation. When $\epsilon_{\mathrm{POD}}=10^{-6}$ and $N_w=2$, $r_{c,1}=r_{c,2}=25$ for both SI-TAR-IG and FGMRES-TAR-IG. The symbol ``$\star$'' indicates failure to converge within the maximum number of transport sweeps. A dash ``--'' indicates that the entry is not applicable.}}
\label{tab:pincell_online}
\begin{subtable}[t]{\textwidth}
\centering
\small
\setlength{\tabcolsep}{4pt}
\vspace{-5pt}
\caption{SI-based methods.}
\label{tab:pincell_SI}
\begin{adjustbox}{max width=\linewidth,center}
\begin{tabular}{@{}lcccccc@{}}
\toprule
 & SI-DSA & SI-DSA-IG & SI-ROMSAD-2,2
 & SI-ROMSAD-IG-2,2 & SI-TAR & SI-TAR-IG \\
\midrule
\multicolumn{7}{c}{(1) $\epsilon_{\mathrm{POD}}=10^{-6}$,
$r_{\mathrm{IG}}=12$, $N_w=2$}\\
\midrule
$\bar{n}_{\textrm{sweep}}$
& 50 $(\star)$ & 35.7 & 33 & 28.9 & 10.9 & 8.6 \\
$\bar{\mathcal{R}}_{\infty}$
& $1.09\times10^{-9}$ & $1.42\times10^{-12}$
& $9.02\times10^{-13}$ & $6.84\times10^{-13}$
& $6.32\times10^{-13}$ & $6.48\times10^{-13}$ \\
$\bar{T}_{\textrm{rel}}$
& 100\% & 67.86\% & 63.17\% & 55.27\%
& 20.49\% & 16.20\% \\
$N_{\mathrm{BE}}\;(N_{\mathrm{BE}}^{\mathrm{reuse}})$
& -- & -- & 75\;(7) & 62\;(6) & 34\;(2) & 32\;(2)\\
\midrule
\multicolumn{7}{c}{(2) $\epsilon_{\mathrm{POD}}=10^{-9}$,
$r_{\mathrm{IG}}=21$, $N_w=1$}\\
\midrule
$\bar{n}_{\textrm{sweep}}$
& 50 $(\star)$ & 13.5 & 13 & 9.3 & 16.8 & 4.8 \\
$\bar{\mathcal{R}}_{\infty}$
& $1.09\times10^{-9}$ & $1.03\times10^{-12}$
& $7.65\times10^{-13}$ & $7.04\times10^{-13}$
& $7.63\times10^{-13}$ & $3.47\times10^{-13}$ \\
$\bar{T}_{\textrm{rel}}$
& 100\% & 23.94\% & 24.93\% & 16.87\%
& 31.43\% & 9.15\% \\
$N_{\mathrm{BE}}\;(N_{\mathrm{BE}}^{\mathrm{reuse}})$
& -- & -- & 32\;(3) & 29\;(3) & 33\;(2) & 25\;(1)\\
\bottomrule
\end{tabular}
\end{adjustbox}
\end{subtable}
\vspace{0.3cm}
\begin{subtable}[t]{\textwidth}
\centering
\small
\setlength{\tabcolsep}{4pt}
\caption{GMRES-based methods.}
\label{tab:pincell_GMRES}
\begin{adjustbox}{max width=\linewidth,center}
\begin{tabular}{@{}lcccc@{}}
\toprule
 & GMRES-DSA & GMRES-DSA-IG
 & FGMRES-TAR & FGMRES-TAR-IG \\
\midrule
\multicolumn{5}{c}{(1) $\epsilon_{\mathrm{POD}}=10^{-6}$,
$r_{\mathrm{IG}}=12$, $N_w=2$}\\
\midrule
$\bar{n}_{\textrm{sweep}}$
& 20.6 & 12.3 & 6.4 & 6.5 \\
$\bar{\mathcal{R}}_{\infty}$
& $5.59\times10^{-13}$ & $6.57\times10^{-13}$
& $8.03\times10^{-13}$ & $6.92\times10^{-13}$ \\
$\bar{T}_{\textrm{rel}}$
& 39.07\% & 22.93\% & 11.71\% & 11.53\% \\
$N_{\mathrm{BE}}\;(N_{\mathrm{BE}}^{\mathrm{reuse}})$
& -- & -- & 30\;(3) & 31\;(4)\\
\midrule
\multicolumn{5}{c}{(2) $\epsilon_{\mathrm{POD}}=10^{-9}$,
$r_{\mathrm{IG}}=21$, $N_w=1$}\\
\midrule
$\bar{n}_{\textrm{sweep}}$
& 20.6 & 6.6 & 7.0 & 4.1 \\
$\bar{\mathcal{R}}_{\infty}$
& $5.59\times10^{-13}$ & $4.84\times10^{-13}$
& $7.11\times10^{-13}$ & $5.88\times10^{-13}$ \\
$\bar{T}_{\textrm{rel}}$
& 39.07\% & 12.31\% & 13.05\% & 7.31\% \\
$N_{\mathrm{BE}}\;(N_{\mathrm{BE}}^{\mathrm{reuse}})$
& -- & -- & 30\;(2) & 29\;(2)\\
\bottomrule
\end{tabular}
\end{adjustbox}
\end{subtable}
}
\end{table}

\textbf{Online performance and \ning{robustness}:} We consider two pairs of hyperparameters for the POD truncation tolerance and the total aware level, $(\epsilon_{\mathrm{POD}}, N_w) = (10^{-6},2)$ and $(10^{-9},1)$. Our main observations from Table~\ref{tab:pincell_online} are as follows. 

\begin{enumerate}
    \item \textit{Robustness with incompletely converged training snapshots.} Standard SI-DSA fails to reach the $10^{-11}$ tolerance within the $50$-sweep limit, with an average residual of $\bar{\mathcal{R}}_{\infty} = 1.09\times10^{-9}$. \ning{Remarkably, despite the limited accuracy of the training snapshots, all solvers using the TAR framework still converge rapidly online, achieving mean \(\ell_\infty\) residuals below $10^{-12}$.}

    \item \textit{Results for SI.} \ning{For $(\epsilon_{\mathrm{POD}}, N_w) = (10^{-6}, 2)$, SI-TAR-IG converges in $8.6$ sweeps on average, reducing sweep counts by factors of $3.84$ and $3.36$ (with wall-clock speedups of $3.90\times$ and $3.41\times$) compared to SI-ROMSAD-2,2 and SI-ROMSAD-IG-2,2, respectively. For $(10^{-9}, 1)$, SI-TAR-IG requires an average of just $4.8$ sweeps, yielding sweep reduction factors of $2.71$ and $1.94$ ($2.72\times$ and $1.84\times$ speedups) over SI-ROMSAD-2,2 and SI-ROMSAD-IG-2,2, respectively.}

    \item \textit{GMRES results.} \ning{For $\epsilon_{\mathrm{POD}}=10^{-6}$ and $10^{-9}$, FGMRES-TAR-IG converges in $6.5$ and $4.1$ sweeps, respectively, reducing sweep counts by factors of $3.17$ and $5.02$ relative to GMRES-DSA, and $1.89$ and $1.61$ relative to GMRES-DSA-IG. The corresponding wall-clock speedups are $3.39\times$ and $5.34\times$ over GMRES-DSA, and $1.99\times$ and $1.68\times$ over GMRES-DSA-IG.} Since the partially consistent DSA solves are inexpensive relative to the high-order $\mathrm{CL}(30,6)$ sweep ($t_{\mathrm{DSA}}/t_{\mathrm{sweep}} = \ning{6.17\%}, \ning{5.55\%},$ and $\ning{3.92\%}$ for SI-DSA, SI-ROMSAD-IG-2,2, and SI-TAR-IG, respectively), the wall-clock speedups closely track the sweep reduction ratios.

    \item \textit{Effect of awareness depth $N_w$ vs.\ POD tolerance $\epsilon_{\mathrm{POD}}$.} \ning{With a zero initial guess, SI-TAR with $(10^{-6}, N_w=2)$ requires significantly fewer sweeps ($10.9$ vs.\ $16.8$) than with $(10^{-9}, N_w=1)$ despite the looser POD tolerance. This suggests that, when the initial error is large, adding a second trajectory-aware correction level is more effective than merely tightening the POD truncation tolerance at the first level.}
\end{enumerate}

\begin{figure}[!ht]
    \centering        \includegraphics[width=0.439\textwidth]{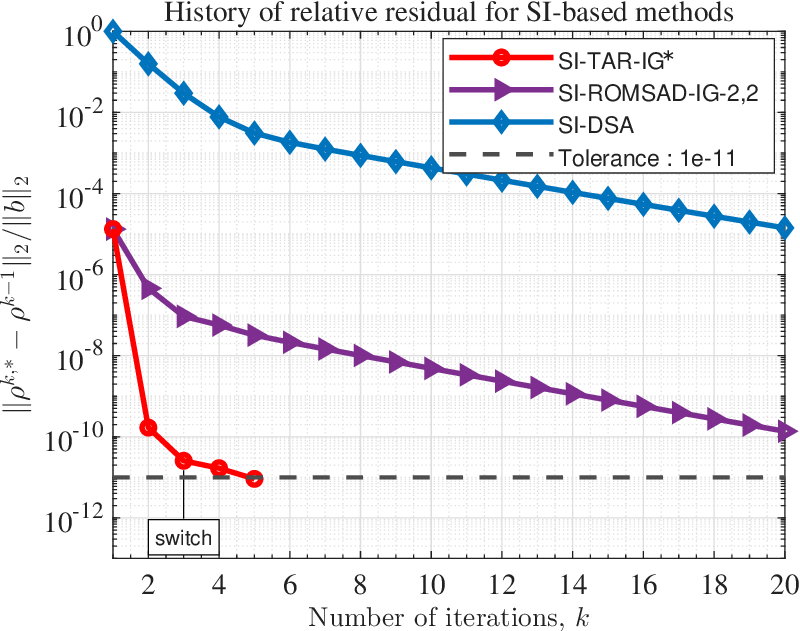}
    \includegraphics[width=0.439\textwidth]{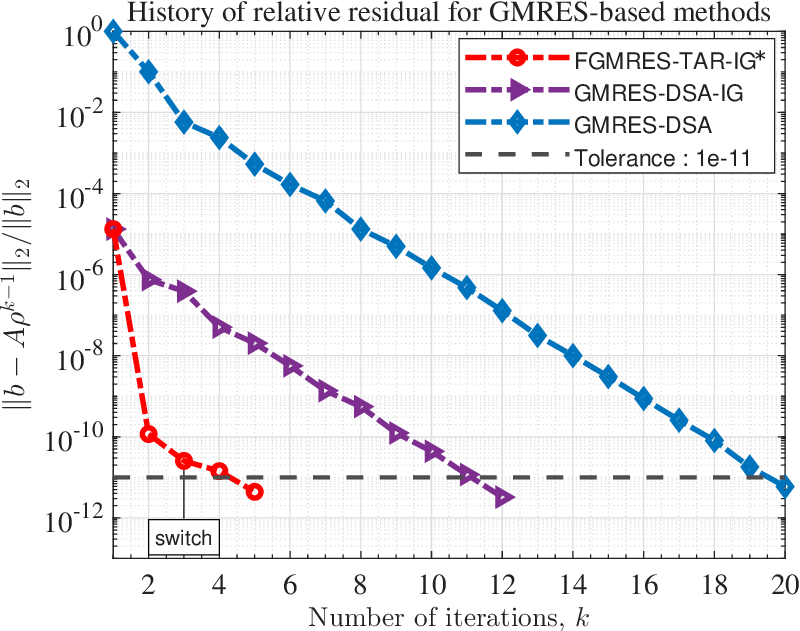}
    \caption{\ning{relative residual histories for the 2D pin-cell problem in Sec.~\ref{sec:pincell} with parameter $(\mu_{a},\mu_{s}) = (0.3550,0.2082)$, where the POD truncation tolerance is $\epsilon_{\mathrm{POD}} = 10^{-6},N_{w}=2$. Left: results for SI-based methods. Right: results for GMRES-based methods. }} 
    \label{fig:pincell}
\end{figure}

\textbf{Convergence history:} Fig.~\ref{fig:pincell} compares the residual histories for a representative test parameter $(\mu_{a},\mu_{s}) = \ning{(0.3550,0.2082)}$. When $\epsilon_{\mathrm{POD}}=10^{-6}$, GMRES-DSA-IG converges faster than \ning{SI-ROMSAD-IG-2,2}, whereas trajectory-aware preconditioning significantly accelerates the convergence of both SI and FGMRES in their initial ROM-preconditioned iterations.

\textbf{Break-even analysis:} \ning{Because SI-DSA fails to converge within $50$ sweeps, we conservatively calculate the break-even counts ($N_{\mathrm{BE}}$ in Table~\ref{tab:pincell_online}) using the computational time of SI-DSA capped at 50 sweeps. Nevertheless, SI-TAR-IG breaks even after only $32$ and $25$ online evaluations for $(\epsilon_{\mathrm{POD}}, N_w) = (10^{-6}, 2)$ and $(10^{-9}, 1)$, respectively. These counts are lower than those of SI-ROMSAD-IG-2,2 ($62$ and $29$) by factors of $1.94$ and $1.16$, and those of SI-ROMSAD-2,2 ($75$ and $32$) by factors of $2.34$ and $1.28$. Both FGMRES-TAR and FGMRES-TAR-IG break even in $29$ to $31$ evaluations across both settings. When high-fidelity training snapshots are reused ($N_{\mathrm{BE}}^{\mathrm{reuse}}$), SI-TAR-IG recovers its offline overhead in just $2$ and $1$ online evaluations, while the FGMRES variants only requires $2$-$4$ solves.}

\begin{table}[!ht]
\centering
\ning{
\begin{threeparttable}
\captionsetup{
    labelformat=simple,
    labelsep=colon
}
\caption{The average relative offline computational cost for the 2D pin-cell problem (Sec.~\ref{sec:pincell}), \ning{where $\epsilon_{\mathrm{POD}}=10^{-9}, N_w=1.$  A dash ``--'' indicates that the entry is not applicable.}}
\label{tab:offline_time_pincell}
\begin{tabular}{lcccc}
\toprule
 & $T_{\mathrm{basis}}$ & $T_{\mathrm{operator}}$ & $\bar{T}_{\mathrm{sweep}}$ & \ning{$\bar{T}_{\mathrm{Arnoldi}}$} \\
\midrule
Initial guess & 0.49\% & 0.36\% &  --  & -- \\
Correction equation (SI-TAR-IG) & 0.74\% & 0.46\% & 1.86\% & -- \\
Correction equation (FGMRES-TAR-IG) & 0.46\% &  0.39\% & 1.65\%\ & 1.61\% \\
\bottomrule\\
\end{tabular}
\end{threeparttable}}
\end{table}
\subsection{Lattice problem \label{sec:lattice}}
\begin{figure}[!ht]
    \centering
    \begin{subfigure}[b]{0.4\textwidth}
        \centering
        \begin{minipage}{6cm}
        \includegraphics[scale=0.5]{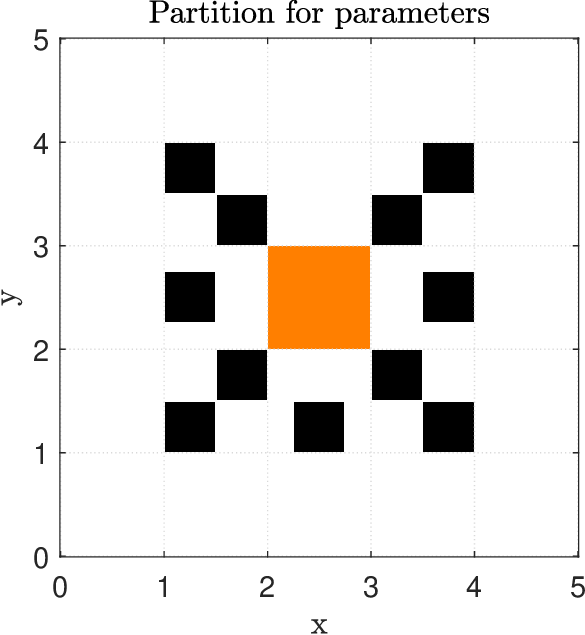}
        \end{minipage}
        \label{fig:lattice-var}
    \end{subfigure}
    \hspace{-0.04\textwidth}
    \begin{subfigure}[b]{0.4\textwidth}
        \centering
        \begin{minipage}{6cm}
        \includegraphics[scale=0.5]{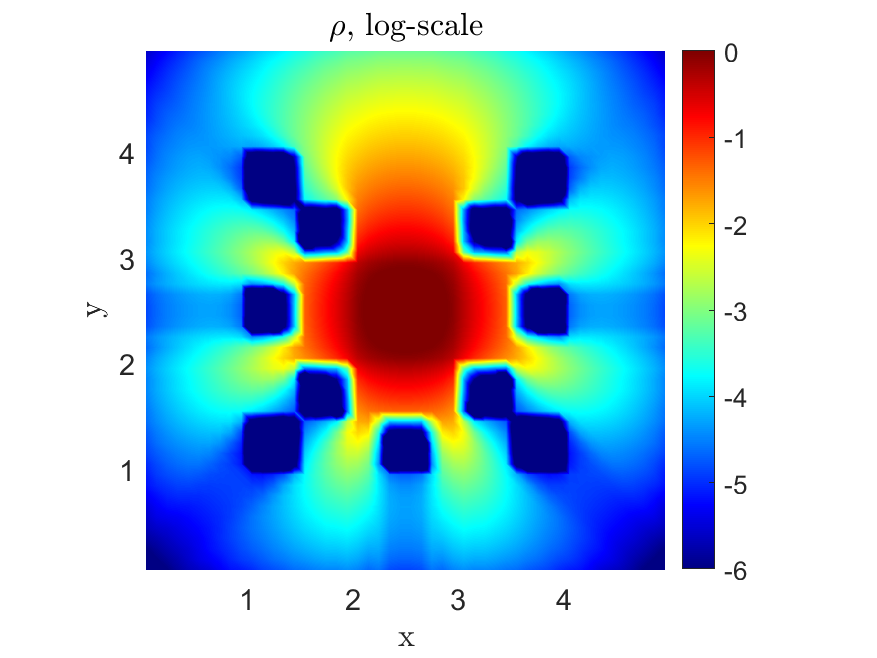}
        \end{minipage}
        \label{fig:lattice-ref}
    \end{subfigure}
    \caption{Left: cross-section setup for the lattice problem in \ning{Sec.~\ref{sec:lattice}}; black regions represent pure absorption zones with $(\sigma_{a},\sigma_{s}) = (\mu_a,0)$, white regions represent pure scattering zones with $(\sigma_{a},\sigma_{s}) = (0,\mu_{s})$, and a unit source is applied in the orange region. Right: reference solution for $(\mu_a,\mu_s) = (97.2353,1.0507)$.} 
    \label{fig:lattice_variable}
\end{figure}

We consider a parametric lattice problem in this section. The computational domain is a square $\mathcal{D}=[0,5]^{2}$ with zero inflow boundary conditions. As shown in Fig.~\ref{fig:lattice_variable}, the black regions represent purely absorbing zones with $(\sigma_{a},\sigma_{s}) = (\mu_a,0)$, while the remaining areas are purely scattering zones with $(\sigma_{a},\sigma_{s}) = (0,\mu_{s})$. The parameter $\mu = (\mu_a,\mu_s)$ controls the strength of absorption and scattering, where $\mu_a \in [95,105]$ and $\mu_s \in [0.5,1.5]$. The source term is nonzero only in the orange region and is defined as
\begin{equation}
G(x, y)=\left\{\begin{array}{ll}
1.0, & \text{if } |x-2.5|<0.5 \text{ and } |y-2.5|<0.5, \\
0, & \text{otherwise.}
\end{array}\right.
\end{equation}
We use a $50\times50$ uniform spatial grid combined with CL(40,6) angular quadrature. The stopping tolerances are set to $\ning{\epsilon_{\textrm{SI}}}=\epsilon_{\mathrm{GMRES}}=10^{-12}$.

\ning{We use this benchmark to investigate the effect of training-set size $N_{\mathrm{train}}$ on the performance of the trajectory-aware framework and the resulting trade-off between offline cost and online speedup. Specifically, we construct uniform $m\times m$ tensor-product training grids over $[95,105]\times[0.5,1.5]$ with $m\in\{7,9,11\}$, yielding $N_{\mathrm{train}}=49,81,$ and $121$ training samples, respectively. For online validation, we randomly select $10$ test parameter pairs $(\mu_a,\mu_s)\in[95,105]\times[0.5,1.5]$.}
Figure~\ref{fig:lattice_variable} (right) shows the reference solution at $(\mu_a,\mu_s)=(97.2353,1.0507)$.

\begin{figure}[!ht]
    \centering
    \includegraphics[width=0.439\textwidth]{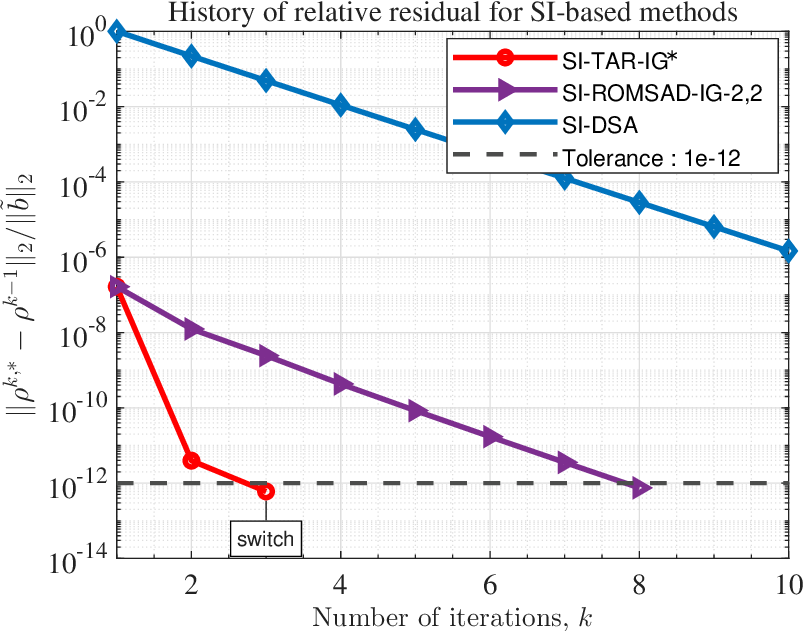}
    \includegraphics[width=0.439\textwidth]{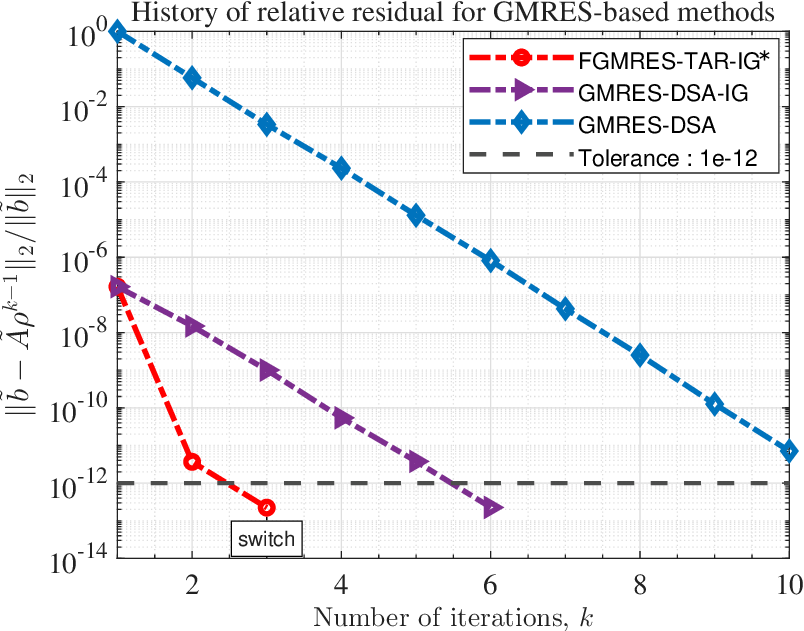}
    \caption{\ning{relative residual histories for the 2D lattice problem in Sec.~\ref{sec:lattice} at $(\mu_a,\mu_s)=(95.3530, 0.6422)$ with $N_{\mathrm{train}}=121$, $\epsilon_{\mathrm{POD}} = 10^{-7}$, and $\ning{N_{w}=2}$. Left: SI-based methods. Right: GMRES-based methods.}} 
    \label{fig:lattice_sample9}
\end{figure}

\ning{\textbf{Baseline comparison ($N_{\mathrm{train}}=121$):}} We first fix $\epsilon_{\mathrm{POD}}=10^{-7}$ and \ning{$N_{w}=2$} on the finest training grid ($N_{\mathrm{train}}=121$, where $r_{\mathrm{IG}}=15$) to benchmark our methods against standard baselines. As shown in Table~\ref{tab:lattice_online} and the residual histories in Fig.~\ref{fig:lattice_sample9}, both trajectory-aware solvers significantly outperform their counterparts:
\begin{enumerate}
    \item \ning{SI-TAR-IG requires $3.4$ transport sweeps on average ($\bar{\mathcal{R}}_{\infty} = 6.15\times10^{-14}$), reducing sweep counts by factors of $6.18$ and $2.76$ and achieving wall-clock speedups of $7.11\times$ and $3.08\times$ over SI-DSA and SI-ROMSAD-IG-2,2, respectively.}
    \item \ning{Similarly, FGMRES-TAR-IG converges in $3.4$ sweeps on average, reducing sweep counts by factors of $3.26$ and $2.06$ (with wall-clock speedups of $3.82\times$ and $2.31\times$) over GMRES-DSA ($11.1$ sweeps) and GMRES-DSA-IG ($7$ sweeps).} 
\end{enumerate}

\begin{table}[!ht]
\centering
\caption{\ning{Baseline comparison for the 2D lattice problem (Sec.~\ref{sec:lattice}) with $N_{\mathrm{train}}=121$, $\epsilon_{\mathrm{POD}}=10^{-7}$, $N_w=2$, and $r_{\mathrm{IG}}=15$. The ratio of DSA time to sweep time ($t_{\mathrm{DSA}}/t_{\mathrm{sweep}}$) is $14.82\%$ for SI-DSA, $9.53\%$ for SI-ROMSAD-IG-2,2, and $0.63\%$ for SI-TAR-IG.}}
\label{tab:lattice_online}
\ning{
\small
\setlength{\tabcolsep}{4pt}
\begin{adjustbox}{max width=\linewidth,center}
\begin{tabular}{@{}lcccccc@{}}
\toprule
 & SI-DSA
 & SI-ROMSAD-IG-2,2
 & SI-TAR-IG
 & GMRES-DSA
 & GMRES-DSA-IG
 & FGMRES-TAR-IG \\
\midrule
$\bar{n}_{\mathrm{sweep}}$
& 21 & 9.4 & 3.4 & 11.1 & 7 & 3.4 \\
$\bar{\mathcal{R}}_{\infty}$
& $2.95\times10^{-13}$
& $8.32\times10^{-14}$
& $6.15\times10^{-14}$
& $1.84\times10^{-13}$
& $1.06\times10^{-13}$
& $1.34\times10^{-13}$ \\
$\bar{T}_{\mathrm{rel}}$
& 100\% & 43.26\% & 14.06\%
& 52.69\% & 31.88\% & 13.79\% \\
\bottomrule
\end{tabular}
\end{adjustbox}
}
\end{table}

\ning{\textbf{Sensitivity to training-set size:} Table~\ref{tab:training-size} compares the performance of all four trajectory-aware variants across $N_{\mathrm{train}} = 49, 81,$ and $121$ with $\epsilon_{\mathrm{POD}}=10^{-7}$ and $N_w=2$. We highlight three key findings:}
\begin{enumerate}
    \item \textit{Offline-online trade-off across $N_{\mathrm{train}}$.} \ning{Increasing the training-set size from $49$ to $121$ enriches the reduced-basis dimensions $(r_{c,1}, r_{c,2})$ and reduces the average online sweeps by factors of $1.93$ and $2.47$ for SI-TAR and SI-TAR-IG, respectively, and by $1.88$ for both FGMRES-TAR and FGMRES-TAR-IG. Conversely, the increased offline cost raises the break-even count $N_{\mathrm{BE}}$ from $78$–$89$ to $166$–$180$ online evaluations (and from $11$–$18$ to $21$–$37$ for $N_{\mathrm{BE}}^{\mathrm{reuse}}$). Even with the smallest training set $(N_{\mathrm{train}}=49)$, the TAR solvers achieve online speedups of $2.5\times$ to $3.6\times$ over SI-DSA and break even within $90$ solves.}

 \pzc{Overall, a larger training set yields greater online acceleration at the expense of a higher break-even count. Hence, moderate training sets are preferred when only a modest number of online solves are needed, whereas larger training sets are preferred when many online solves are required.}

    \item \textit{Contrasting sensitivity of SI and FGMRES to initial guesses.} \ning{As $N_{\mathrm{train}}$ increases, the advantage of SI-TAR-IG over SI-TAR with a zero initial guess becomes more pronounced (reducing sweeps from $4.4$ to $3.4$ at $N_{\mathrm{train}}=121$). In contrast, FGMRES-TAR slightly outperforms FGMRES-TAR-IG in wall-clock time across all training-set sizes. Although FGMRES-TAR-IG requires fewer Krylov iterations ($1.4$ vs. $2.3$ at $N_{\mathrm{train}}=121$), this reduction is offset by the additional transport sweep needed to compute the nonzero initial residual.}

    \item \textit{Offline scaling at large $N_{\mathrm{train}}$.} \ning{Table~\ref{tab:offline_time_lattice} breaks down the relative offline overhead at $N_{\mathrm{train}}=121$. Because SVD is performed on larger snapshot matrices, basis construction ($T_{\mathrm{basis}}$) accounts for $24.13\%$ to $51.72\%$ of a single SI-DSA solve, while projected-operator construction requires at most $9.08\%$. With two aware levels, the additional transport sweeps ($\bar{T}_{\mathrm{sweep}}$) cost $7.59\%$–$9.38\%$, while the Arnoldi updates ($\bar{T}_{\mathrm{Arnoldi}}$) for FGMRES-TAR-IG account for $8.19\%$ of a single SI-DSA solve. Although these costs increase compared with the single-level setting, they remain manageable relative to one linear solve.}
\end{enumerate}

\begin{table}[!ht]
\centering
\ning{
\begin{threeparttable}
\captionsetup{
    labelformat=simple,
    labelsep=colon
}
\caption{Sensitivity to training-set size $N_{\mathrm{train}}$
for the 2D lattice problem (Sec.~\ref{sec:lattice}) with
$\epsilon_{\mathrm{POD}}=10^{-7}$ and $N_w=2$.
Here, $r_{c,l}$ denotes the reduced-basis dimension at the
$l$-th aware level.}
\label{tab:training-size}
\setlength{\tabcolsep}{4pt}
\renewcommand{\arraystretch}{1.0}
\begin{tabular}{lccccccc}
\toprule
Method & $N_{\mathrm{train}}$
& $r_{c,1}$ & $r_{c,2}$
& $\bar{n}_{\mathrm{sweep}}$
& $\bar{\mathcal{R}}_{\infty}$
& $\bar{T}_{\mathrm{rel}}$
& $N_{\mathrm{BE}}\;(N_{\mathrm{BE}}^{\mathrm{reuse}})$ \\
\midrule
\multirow{3}{*}{SI-TAR}
& 49  & 15 & 38 & 8.5
& $1.08\times10^{-13}$ & $39.71\%$ & $89\;(11)$ \\
& 81  & 15 & 43 & 6
& $7.95\times10^{-14}$ & $26.19\%$ & $123\;(14)$ \\
& 121 & 15 & 47 & 4.4
& $5.12\times10^{-14}$ & $18.48\%$ & $172\;(21)$ \\
\midrule
\multirow{3}{*}{SI-TAR-IG}
& 49  & 37 & 49 & 8.4
& $7.73\times10^{-14}$ & $39.27\%$ & $89\;(12)$ \\
& 81  & 43 & 81 & 5.4
& $6.31\times10^{-14}$ & $23.44\%$ & $120\;(15)$ \\
& 121 & 47 & 121 & 3.4
& $6.15\times10^{-14}$ & $14.06\%$ & $166\;(23)$ \\
\midrule
\multirow{3}{*}{FGMRES-TAR}
& 49  & 15 & 39 & $6.2\;(5.2+1)$
& $9.66\times10^{-14}$ & $27.52\%$ & $78\;(13)$ \\
& 81  & 15 & 45 & $4.2\;(3.2+1)$
& $1.76\times10^{-13}$ & $17.88\%$ & $116\;(18)$ \\
& 121 & 15 & 50 & $3.3\;(2.3+1)$
& $1.19\times10^{-13}$ & $13.52\%$ & $172\;(29)$ \\
\midrule
\multirow{3}{*}{FGMRES-TAR-IG}
& 49  & 38 & 49 & $6.4\;(4.4+2)$
& $2.67\times10^{-13}$ & $27.91\%$ & $83\;(18)$ \\
& 81  & 43 & 81 & $4.9\;(2.9+2)$
& $1.59\times10^{-13}$ & $20.86\%$ & $127\;(25)$ \\
& 121 & 47 & 121 & $3.4\;(1.4+2)$
& $1.34\times10^{-13}$ & $13.79\%$ & $180\;(37)$ \\
\bottomrule
\end{tabular}
\end{threeparttable}
}
\end{table}

\begin{table}[!ht]
\centering
\ning{
\begin{threeparttable}
\captionsetup{
    labelformat=simple,
    labelsep=colon
}
\caption{Average relative offline computational cost for the 2D lattice problem (Sec.~\ref{sec:lattice}) at $N_{\mathrm{train}}=121$ with $\epsilon_{\mathrm{POD}}=10^{-7}$ and $N_w=2$.  A dash ``--'' indicates that the entry is not applicable.}
\label{tab:offline_time_lattice}
\begin{tabular}{lcccc}
\toprule
 & $T_{\mathrm{basis}}$ & $T_{\mathrm{operator}}$ & $\bar{T}_{\mathrm{sweep}}$  & \ning{$\bar{T}_{\mathrm{Arnoldi}}$} \\
\midrule
Initial guess & 24.13\% & 0.62\% &  -- & -- \\
Correction equation (SI-TAR-IG) & 51.43\% & 8.87\% & 9.38\% & --\\
Correction equation (FGMRES-TAR-IG) & 51.72\% & 9.08\% & 7.59\% & 8.19\% \\
\bottomrule
\end{tabular}
\end{threeparttable}}
\end{table}

\section{Conclusions\label{sec:conclusion}}
Recognizing the efficiency reduction in ROMSAD caused by the mismatch between offline and online residual trajectories, we propose \ning{an SA-derived, trajectory-aware ROM
preconditioning framework tailored to parametric RTEs:}
\begin{enumerate}
    \item For SI, we construct trajectory-aware ROMs by sequentially building a reduced-order space with ROM-based corrections from earlier iterations. This approach ensures consistency between the offline and online residual trajectories.
    \item We further extend our framework to FGMRES by \ning{reformulating the ideal kinetic correction equation and constructing iteration-dependent ROM preconditioners using the Arnoldi iteration.}
\end{enumerate}

Numerical results validate the effectiveness of the proposed trajectory-aware framework with the following main observations:
\begin{enumerate}
    \item Leveraging low-rank structures across parameters, the online stage of our method achieves a significant acceleration over DSA. 
    \item With only modest additional offline cost compared to ROMSAD, the trajectory-aware framework substantially \ning{mitigates the offline-online mismatch} and \pzc{greatly improves} efficiency in the online stage, especially under loose POD truncation tolerances. \ning{The low break-even numbers further show an efficient offline-online trade-off, particularly when the high-fidelity training snapshots are already available.}
\end{enumerate}

Several directions remain for future research.  First, to the best of our knowledge, a complete theoretical understanding of the relation between the dimension of ROMs and the convergence rate of iterative solvers with ROM-based preconditioners is still lacking and highly desirable, even for parametric elliptic equations. Second, we aim to extend our trajectory-aware framework to more complex applications, including non-affine or time-dependent problems, nonlinear thermal radiation, and problems with multi-group energy and anisotropic scattering. In the longer term, we plan to integrate this framework as a building block for uncertainty quantification, design optimization, and inverse problems. 

\section*{Acknowledgment}
This work was partially supported by the Hong Kong Research Grants Council grants Early Career Scheme 26302724 and General Research Fund 16306825. Besides the funding agency, the authors would also like to thank Dr. Ningxin Liu for her kind proofreading and suggestions.

\section*{CRediT authorship contribution statement}

{\bf Ning Tang:} Writing – original draft, Writing – review \& editing, Visualization, Validation, Software, Methodology, Data curation, Conceptualization.

{\bf Zhichao Peng:} Writing – review \& editing, Visualization, Validation, Software, Methodology, Data curation, Conceptualization, Supervision, Funding acquisition.

\section*{Declaration of generative AI and AI-assisted technologies in the writing process}
During the preparation of this work, the authors used AI to check for grammar errors, improve readability \ning{and assist with code development and debugging.} After using this tool/service, the authors reviewed and edited the content as needed and take full responsibility for the content of the publication.

\appendix
\section{Reformulation of the ideal correction equation\label{sec:equation_equivalence}}

As discussed in Sec.~\ref{sec:ideal_correction_krylov}, the ideal correction equation for FGMRES is 
\begin{equation}
(\BD_j+\BSigma_t)\delta\bff_j^{(l)}=\BSigma_s\delta\brho^{(l)}+\BSigma_s \bq^{(l)},\quad \delta\brho^{(l)}=\sum_{j=1}^{N_{\Bupsilon}}\omega_j\delta\bff^{(l)}_j,\quad j=1,\dots,N_{\Bupsilon}.  
\label{eq:ideal-kinetic-correction-fgmres}
\end{equation} 
To avoid directly solving this equation when building ROMs, we follow the approach in our previous work \cite{peng2025flexible}. Define $\beeta^{(l)}$ as the solution to 
\begin{equation}
    \WBA\beeta^{(l)}=(\BI-\mathbf{K}\BSigma_s)\beeta^{(l)}=\bq^{(l)}.\label{eq:xi-def}
    \end{equation}
Substituting \eqref{eq:xi-def} into \eqref{eq:ideal-kinetic-correction-fgmres}, we obtain
\begin{align}
     (\BD_j+\BSigma_t)\delta\bff_j^{(l)}&=
     \BSigma_s\delta\brho^{(l)}+\BSigma_s(\BI-\BK\BSigma_s)\beeta^{(l)},\\
     \delta\bff_j^{(l)}&=(\BD_j+\BSigma_t)^{-1}\left[\BSigma_s\delta\brho^{(l)}+\BSigma_s(\BI-\BK\BSigma_s)\beeta^{(l)}\right].
     \label{eq:density_correction_appendix}
\end{align}
Then, the macroscopic density correction, $\delta\brho^{(l)}$, can be computed as:
\begin{align}
\delta\brho^{(l)}&=\big(\sum_{j=1}^{N_{\Bupsilon}}\omega_j (\BD_j+\BSigma_t)^{-1}\big)\left[\BSigma_s\delta\brho^{(l)}+\BSigma_s(\BI-\BK\BSigma_s)\beeta^{(l)}\right].
\end{align}
Based on the definition of $\mathbf{K}=\sum_{j=1}^{N_{\Bupsilon}}\omega_j (\BD_j+\BSigma_t)^{-1}$ in \eqref{eq:matrix_vec_rewrite_rho}, we obtain
\begin{subequations}
\begin{align}
   \delta\brho^{(l)}&= \BK\BSigma_s\delta\brho^{(l)}+\BK\BSigma_s(\BI-\BK\BSigma_s)\beeta^{(l)}\\
   (\BI-\BK\BSigma_s)\delta\brho^{(l)}&=\BK\BSigma_s(\BI-\BK\BSigma_s)\beeta^{(l)}=(\BI-\BK\BSigma_s)\BK\BSigma_s\beeta^{(l)}.
 \end{align}
\end{subequations}
Hence, 
\begin{equation}
\delta\brho^{(l)}=\BK\BSigma_s\beeta^{(l)}.
\label{eq:density-relation}
\end{equation}
Substituting \eqref{eq:density-relation} into \eqref{eq:density_correction_appendix}, we have
\begin{subequations}
\begin{align}
(\BD_j+\BSigma_t)\delta\bff_j^{(l)}&=\BSigma_s\BK\BSigma_s\beeta^{(l)}+\BSigma_s(\BI-\BK\BSigma_s)\beeta^{(l)},\\
(\BD_j+\BSigma_t) \delta\bff_j^{(l)} &=\BSigma_s\beeta^{(l)},\;\; j=1,\dots,N_{\Bupsilon}, 
    \label{eq:equivalence_1_appendix}
\end{align}
\end{subequations}
where $\beeta^{(l)}$ satisfies
$ (\BI-\mathbf{K}\BSigma_{s})\beeta^{(l)} = \bq^{(l)}.
    \label{eq:equivalence_2_appendix}
$

\section{Equivalence between SA and left preconditioning\label{sec:left_preconditioner}}
After applying the SA correction given by
\begin{equation}
\mathbf{C}\delta\brho^{(l)} = \BSigma_s(\brho^{(l,*)}-\brho^{(l-1)})=\BSigma_s\br^{(l-1)},
\label{eq:delta_rho_c}
\end{equation}the SI update becomes:
\begin{equation}
\begin{split}
    \brho^{(l)} = &\brho^{(l,*)} + \delta\brho^{(l)}\\
    =& \brho^{(l-1)} + \boldsymbol{r}^{(l-1)} +\mathbf{C}^{-1}\BSigma_s\boldsymbol{r}^{(l-1)} \\
    =& \brho^{(l-1)} + (\mathbf{I}+\mathbf{C}^{-1}\BSigma_s)\boldsymbol{r}^{(l-1)}.
\end{split}
\end{equation}
Here, $(\mathbf{I}+\mathbf{C}^{-1}\BSigma_s)\boldsymbol{r}^{(l-1)}$ is the residual of the left-preconditioned system
$(\mathbf{I}+\mathbf{C}^{-1}\BSigma_{s})(\mathbf{I} - \mathbf{K}\BSigma_s) \brho = (\mathbf{I}+\mathbf{C}^{-1}\BSigma_{s})\widetilde{\mathbf{b}}.$
Hence, SA is equivalent to introducing a left preconditioner $\BM^{-1}=\BI+\BC^{-1}\BSigma_s$.


\section{Condition numbers for reduced systems \label{sec:condition_number}}
\ning{For each aware level $l$, we define the mean condition number of the Galerkin reduced matrices over the test set
$\mathscr{P}_{\mathrm{test}}$ as
\begin{equation}
\bar{\kappa}_{c,l}
= \frac{1}{N_{\mathrm{test}}}
\sum_{\mu\in\mathscr{P}_{\mathrm{test}}}
\kappa_2\left(
(\BU_{r_l}^{(l)})^T \BA_\mu \BU_
{r_l}^{(l)}
\right),
\end{equation}
where $N_{\mathrm{test}}=|\mathscr{P}_{\mathrm{test}}|$ and $\kappa_2$ denotes the condition number in the Euclidean norm.}
\ning{Table~\ref{tab:condition_2d} reports the mean and maximum condition numbers for the 2D problems considered in this work. The largest values occur in the variable-scattering problem, and all are at  most of order $10^3$. In double precision, their products with machine epsilon are of order $10^{-13}$ or smaller, suggesting that round-off errors are unlikely to be significantly amplified when solving the reduced systems with backward-stable solvers.}
\begin{table}[!ht]
\centering
\small
\ning{
\caption{Mean condition numbers of the reduced systems for the 2D problems, with maxima in parentheses.}
\label{tab:condition_2d}
\begingroup
\renewcommand{\arraystretch}{1.2}
\setlength{\tabcolsep}{5pt}
\begin{tabularx}{.9\linewidth}{
@{} c *{4}{>{\centering\arraybackslash}X} @{}
}
\toprule
& SI-TAR & SI-TAR-IG & FGMRES-TAR & FGMRES-TAR-IG \\
\midrule
\multicolumn{5}{c}{(1) Variable-scattering (Sec.~\ref{sec:scattering}): $\epsilon_{\mathrm{POD}}=10^{-7},\;\;
 N_w=1$} \\
\midrule
\addlinespace[3pt]
$\bar{\kappa}_{c,1}$
& 892.16 (1609.02)
& 1298.02 (2355.24)
& 892.00 (1608.71)
& 1325.07 (2419.49) \\
\midrule
\multicolumn{5}{c}{(2) Pin-cell (Sec.~\ref{sec:pincell} ): $\epsilon_{\mathrm{POD}}=10^{-9},\;\;
 N_w=1$} \\
\midrule
\addlinespace[3pt]
$\bar{\kappa}_{c,1}$
& 51.74 (115.21)
& 67.34 (148.15)
& 51.75 (115.23)
& 47.54 (104.64) \\
\midrule
\multicolumn{5}{c}{(3) Lattice (Sec.~\ref{sec:lattice}): $\epsilon_{\mathrm{POD}}=10^{-7},\;\;
 N_w=2$} \\
\midrule
\addlinespace[3pt]
$\bar{\kappa}_{c,1}$
& 171.55 (188.29)
& 162.23 (180.81)
& 173.25 (190.13)
& 161.38 (179.83) \\
$\bar{\kappa}_{c,2}$
& 169.08 (187.62)
& 64.56 (77.41)
& 138.11 (151.43)
& 111.39 (124.51) \\
\bottomrule
\end{tabularx}
\endgroup}
\end{table}

\bibliographystyle{elsarticle-num} 
\bibliography{ref}
\end{document}